\documentclass[a4paper,12pt]{amsart}
\usepackage[english]{babel}
\usepackage{latexsym}
\usepackage{amsmath}
\usepackage{amssymb}	
\usepackage{amsthm}
\usepackage{xcolor}

\usepackage[textwidth=400pt, textheight=620pt]{geometry}
\usepackage{enumitem}

\newtheorem{thm}{Theorem}[section]
\newtheorem{cor}[thm]{Corollary}

\newtheorem{prop}[thm]{Proposition}
\newtheorem{claim}[thm]{Claim}

\theoremstyle{plain}

\newtheoremstyle{def}{3mm}{2mm}{}{0pt}{\bfseries}{}{ }
{\thmname{#1}\thmnumber{ #2 }\thmnote{(#3)} 
}
\theoremstyle{def}  
\newtheorem{defi}[thm]{Definition}
\newtheorem{remark}[thm]{Remark}
\newtheorem{lemma}[thm]{Lemma}

\usepackage{hyperref}

\newcommand{\n}[2]{\Vert #1\Vert_{#2}}

\renewcommand\t{\theta}
\newcommand\di{{\delta^\prime(\t)}}
\newcommand\dn{{\delta^{(n)}(\t)}}

\renewcommand\l{{\lambda}}

\newcommand\fn{f^{(n)}(\t)}
\newcommand\fm{f^{(m)}(\t)}

\newcommand\bii{\bar B_{\di}}
\renewcommand\Re{\operatorname{Re}}
\renewcommand\Im{\operatorname{Im}}
\newcommand\BC{\operatorname{BC}}
\newcommand\C{\operatorname{C}}

\newcommand\bC{\mathbb C}
\newcommand\subC[1]{{(\bC;\lambda\bC)_\theta}^{(#1)}}
\newcommand\cF{\mathcal F}

\begin{document}

\title[New results on Fourier Multipliers on $L^p$]{New results on Fourier Multipliers on $L^p$: \\ a perspective through Unimodular Symbols} 
\author {Mar\'{\i}a J. Carro\and  Alberto Salguero-Alarc\'on} 

\address{Mar\'{\i}a J. Carro.  Department of  Analysis and Applied Mathematics,              
Universidad Complutense de Madrid, Plaza de las Ciencias 3,
28040 Madrid, Spain. Instituto de Ciencias Matem\'aticas ICMAT, Madrid, Spain}
\email{mjcarro@ucm.es}

\address{Alberto Salguero-Alarc\'on.  Department of  Analysis and Applied Mathematics,              
Universidad Complutense de Madrid, Plaza de las Ciencias 3,
28040 Madrid, Spain. Instituto de Ciencias Matem\'aticas ICMAT, Madrid, Spain}
\email{albsalgu@ucm.es}

\date{\today}

 


\thanks{The first author is partially supported by grants 
PID2024-155917NB-I00 funded by MCIN/AEI/10.13039/501100011033,  CEX2019-000904-S funded by MCIN/AEI/ 10.13039/501100011033 and Grupo UCM-970966 (Spain). The second author is partially supported by grants PID2024-162214NB-I00 and PID2023-146505NB-C21 funded by MICIU/AEI/ 10.13039/501100011033.}

\begin{abstract}
The paper focuses on the behaviour of unimodular Fourier multipliers with exponential growth in the context of weighted $L^p$-spaces. Our main result shows that much of the general theory of multipliers is approachable through the theory of unimodular multipliers. Indeed, we show that a bounded  measurable function $m$ is a multiplier on $L^p$ for $1\leq p < \infty$ provided that $e^{itm}$ is a multiplier on $L^p$ $(1\leq p<\infty)$ and its multiplier norm admits an exponential bound of the form $e^{c|t|^s}$ for suitable $c>0$ and $0<s<1$. We then apply this principle to obtain new results related to the boundedness of homogeneous rough operators, singular operators along curves and oscillatory integrals. A key ingredient in our study is an extension of the classical Stein's theorem on analytic families of operators that studies the behaviour of the derivative operator when $\theta \to 0$. 	
\end{abstract}

\maketitle

\pagestyle{headings}\pagenumbering{arabic}\thispagestyle{plain}

\tableofcontents

\section{Introduction} 

\par Given a bounded measurable function $m$  on $\mathbb R^n$, consider, initially on test functions,  the linear operator $T_m$ given by 
\begin{equation} 
	\label{eq:Tm} T_mf(x)= \int_{\mathbb R^n} m(\xi) \hat f(\xi) e^{ix\xi} d\xi,
\end{equation}
where $\hat f$ is the Fourier transform of $f$. The function $m$ is called the \emph{symbol} or \emph{multiplier}.
In fact, if $T_m$ can be extended to a bounded operator on $L^p(\mathbb R^n)$ for a given $1\leq p< +\infty$, we say that $m$ is a \emph{Fourier multiplier} on $L^p(\mathbb R^n)$ and we write  $m\in \mathcal M_p(\mathbb R^n)$. It is a classical fact that $\mathcal M_p(\mathbb R^n)$ is a Banach space (in fact, a Banach algebra under the pointwise product) with the norm $\|m\|_{\mathcal M_p} = \|T_m\|_{L^p \to L^p}$. The characterization of $\mathcal M_p(\mathbb R^n)$ is easy in the cases $p=1$ and $p=2$, and a completely open problem for the other values of $p$. In particular, the Plancherel theorem asserts that $\mathcal  M_2(\mathbb R^n)=L^\infty(\mathbb R^n)$, and it is known that $\mathcal M_1(\mathbb R^n)$ coincides with the algebra $\mathcal A(\mathbb R^n)$ of functions which are the Fourier transform of a finite measure. Although $m$ can be taken complex valued, we can always consider its real and its imaginary part,  and therefore in this paper we shall always consider $m$ to be real-valued. Our first result in this context is the following:

\begin{thm}
	\label{sector1} 
	Let $1\leq p<\infty$ and $0<\varphi<\pi/2$. The following assertions are equivalent:
	\begin{enumerate}
		\item[(i)] $m\in  \mathcal M_p$.
		\item[(ii)] There exists $c>0$ such that for every $z$ with $|{\arg(z)}| = \varphi$, we have
		\[ e^{\pm zm}\in  \mathcal M_p, \quad \|e^{\pm z m} \|_{ \mathcal M_p}\lesssim e^{c|z|}.\]
	\end{enumerate}
\end{thm}

We observe that \emph{(i)} trivially implies \emph{(ii)} even for \emph{all} $z\in \mathbb C$. However, the converse is more subtle and stems from the following relation with the theory of \emph{unimodular Fourier multipliers} (that is, multipliers of the form $m(\xi)= e^{i h(\xi)}$ with $h$ a real measurable function):

\begin{thm} 	\label{boundedm1}
	Let $m$ be a bounded function and $1 \le p<\infty$. If there exist $s\in (0,1)$ and $c>0$ such that 
	$$
	\|e^{ itm} \|_{ \mathcal M_p}\lesssim e^{c |t|^{s}}, \quad \forall t\in\mathbb R\qquad\implies\qquad m\in  \mathcal M_p.
	$$
\end{thm}

To frame our previous theorem, let us give a short overview on unimodular Fourier multipliers. An interesting starting point (chronologically speaking) for this theory is the Beurling-Helson theorem \cite{bb:hh}, which asserts that for a real function $\varphi$, 
\begin{itemize}
	\item If $
	||e^{in\varphi}||_{\mathcal A(\mathbb R)}= O(1)$ for every $n\in\mathbb Z$, then $\varphi$ coincides modulo $2\pi$ with some linear  function. 
	
	\item If $||e^{i\lambda \varphi}||_{\mathcal A(\mathbb R)}= O(1)$ for every $\lambda\in\mathbb R$, then $\varphi$ coincides  with some linear  function.  
\end{itemize}
The Beurling-Helson theorem enabled to solve Levy's problem on the description of the endomorphisms of the algebra $\mathcal A(\mathbb R^n)$, proving that only trivial changes of variable are admissible.
This motivated to study the class of functions $\varphi$ such  that
\begin{equation}\label{lloo}
	||e^{i\lambda \varphi}||_{\mathcal M_p}= O(1), \qquad \lambda\in\mathbb R. 
\end{equation}
H\"ormander \cite{ho:ho} proved that any $\varphi\in C^2$ satisfying \eqref{lloo} is linear and conjectured that $C^1$ suffices, which was later confirmed in \cite{ll:oo2}. Hence it is natural to wonder what can be said for a measurable $\varphi$. Several interesting results have been given in this direction, see \cite{ll:oo} for a comprehensive account. In fact, we provide here a new result in connection with the main theorem  in \cite{ll:oo} (see Theorem \ref{edm} and Corollary \ref{edm2} in Section \ref{fm}). 

\par Unimodular Fourier multipliers also arise in connection to the Cauchy problem for dispersive equations. For example, the solution $u(t, x)$ of 
\begin{equation*}
	\begin{cases}
		i\partial_t u+(-\Delta)^{\alpha/2}=0 \\ 
		u(0, x)= u_0(x),
	\end{cases}
\end{equation*}
can be expressed as $u(t,x)= T_{m}u_0(x)$ with $m(\xi)=e^{it|\xi|^\alpha}$. Therefore, the behaviour of the norm of this ``multiplier'' when $t\to +\infty$ has its own interest in this theory. Let us mention that the cases $\alpha = 1, 2$ and $3$ possess special interest since they are related to the wave equation, the Schr\"odinger equation and the Airy equation, respectively.

In light of the above considerations, if a function $m$ satisfies that $e^{im}\in \mathcal M_p$ for some $p\neq 2$ and $m$ is not linear, then necessarily $||e^{i\lambda m}||_{ \mathcal M_p} \to +\infty$ when   $|\lambda|\to +\infty$. 
On the other hand, it is immediate to prove that
\begin{equation} \label{eq:converse}
	m\in  \mathcal M_p \qquad\implies\qquad e^{\pm itm} \in  \mathcal M_p, \quad ||e^{\pm itm} ||_{M_p}\lesssim e^{|t|\,  ||m||_{ \mathcal M_p}}, \ t\in\mathbb R.
\end{equation}
Therefore, our Theorem \ref{boundedm1} can be regarded as a slightly weaker version of the converse implication. We conjecture, however, that the converse of \eqref{eq:converse} (assuming boundedness of $m$) is false; equivalently, Theorem \ref{boundedm1} does not hold for $s=1$. This suspicion is motivated by the fact that Theorem \ref{boundedm1} is obtained as a consequence of Theorem \ref{thm:main-2}, and the latter does not hold for the parameter $s=1$. Nonetheless, our counterexample --see Remark \ref{rem:optimal}-- does not correspond with the above situation.

\par Our technique is strongly based on an extension of Stein's theorem on \emph{analytic families of operators} (AFO, for short) $\{T_z\}_z$ \cite{st:st, CJ} which analyses the behaviour of the \emph{derivative} operator $T_z'$. To motivate the main ideas, let us start with an easy observation. Assume that $(T_z)_{z}$ is defined in the vertical unit strip $\mathbb S = \{ z: 0<\Re z<1\}$. Then, under suitable hypothesis on the boundedness of $T_{it}: L^{p_0} \to L^{q_0}$ and $T_{1+it}: L^{p_1} \to L^{q_1}$, Stein's interpolation theorem applied to $\{T_z\}_z$ concludes that
$$
T_{1/2}: L^{p(1/2)}\longrightarrow L^{q(1/2)}, 
$$
where 
$$
\frac 1{p(1/2)}= \frac 1{2p_0}+\frac 1{2p_1},  \qquad \frac 1{q(1/2)}= \frac 1{2q_0}+\frac 1{2q_1}.
$$
Of course, if $(T_z)_z$ is instead defined in $\mathbb S_r=\{z\in\mathbb C: 0< \Re z <r\}$, for some $r>1$, and now $T_{it}$ and $T_{r+it}$ satisfy the corresponding boundedness conditions, then an application of Stein's theorem to the family  $R_z=T_{rz}$ yields
$$
T_{1/2}=R_{1/(2r)}: L^{p(1/(2r))}\longrightarrow T^{q(1/(2r))}, 
$$
where now
$$
\frac 1{p(1/(2r))}= \frac {1-\frac 1{2r}}{p_0}+\frac {\frac 1{2r}}{p_1},  \qquad \frac 1{q(1/(2r))}= \frac {1-\frac 1{2r}}{q_0}+\frac {\frac 1{2r}}{q_1}. 
$$
From here, it seems natural to conjecture  that, under hopefully mild conditions,  if the family of operators is defined in
$$
\Omega=\{z\in\mathbb C: \Re z>0\}, 
$$
we shall be able to let $r\to +\infty$ and conclude that $T_{1/2}$ inherits  the same boundedness property as $T_0$; that is, 
$$
T_{1/2}:L^{p_0}\longrightarrow L^{q_0}.
$$
This idea has been materialized a number of times to obtain boundedness conditions for homogeneous rough operators  \cite{F}, Hilbert transforms along curves \cite{nrw1, NRW} and oscillatory integrals \cite{ps0:ps0}  on $L^p$ spaces, as well as to obtain extensions of classical ergodic theorems \cite{ns:ns, ns2:ns2}. However, in all of the cited cases, the absence of quantitative bounds for the norms of either $T_{it}$ or $T_{\alpha +it}$ must be counterweighted by proving (or assuming)  boundedness of the corresponding collection of operators on $L^p$ for \emph{every} $p>1$. 

\par The original approach in this paper is that such hypotheses concerning boundedness can be substantially weakened if a suitable bound with respect to $t\in \mathbb R$ is known for $T_{it}$ and $T_{\alpha +it}$ for $\alpha\gg 0$. Indeed, our main technical result essentially states that, if $(T_z)_{z\in \Omega}$ and 	
$T_{it}: L^{p_0} \to L^{q_0}$ and $T_{\alpha + it}: L^{p_1} \to L^{q_1}$ are bounded for every $\alpha \gg0$, with the bounds being of exponential order on $|t|^s$ for some $0<s<1$, then not only $T_{1/2}$ but also --and this is the key point for the applications-- the \emph{derivative} $T'_{1/2}$ is bounded from $L^{p_0}$ to $L^{q_0}$. Actually, the result is also true for the succesive derivatives $T^{n)}_\theta$. However, we shall restrict to the case $n=1$ for simplicity and because all our examples are included in this case. 

\par Most of our applications bear close relation to the theory of unimodular multipliers. Before passing to a description of our main results, let us mention however that our refinement of Stein's theorem  can be applied to much more general operators than Fourier multipliers and to more general spaces than the $L^p$-spaces. In the applications, we shall pay special attention to the case of weighted $L_p$-spaces. 

\subsection{Homogeneous rough operators}

For a given $n\in \mathbb N$, $\Sigma^{n-1}$ denotes the unit sphere of $\mathbb R^n$. The homogeneous rough operator is defined by
\begin{equation} \label{eq:rough}
	S_\Omega f(x)=p.v.  \int_{\mathbb R^n} \frac{\Omega(y')}{|y|^n} f(x-y) dy,
\end{equation}
where $y' = y/|y|$ and $\Omega\in L^1(\Sigma^{n-1})$ satisfies the \emph{cancellation property}; that is,  $\int_{\Sigma^{n-1}} \Omega(y')d\sigma(y')=0$.  Boundedness properties of this operator both in $L^p$ and $L^p(w)$ when $w$ belongs to a Muckhenhoupt class $A_p$ for some $1\leq p < +\infty$ (see Section \ref{sec:weighted} for the definition) have a long history.  In the $L^p$ case, the Calder\'on-Zygmund method of rotations \cite{cz:cz}  proves that if $\Omega \in L \log  L(\Sigma^{n-1})$, then, for every $p>1$
\begin{equation}\label{tomega1}
	S_\Omega: L^p \longrightarrow L^p. 
\end{equation}
By contrast, although considerable effort has been devoted to the case of weighted $L^p$ spaces --see e.g. \cite{ccdo:ccdo, d:d, hrt:hrt} and the references therein--, the theory is still far from complete. In particular, it is an open question to characterize the weights $v$ such that, for every $\Omega\in L^2(\Sigma^{n-1})$ with the cancellation property, $S_\Omega$ is bounded on $L^2(v)$. As far as we know, the best result in this case is that
\begin{equation*}
	S_\Omega: L^2(v) \longrightarrow  L^2(v), \qquad \forall v\in A_1.
\end{equation*}

\par The present paper strengthens a connection between the boundedness of \eqref{eq:rough} and the asymptotic behaviour of $L^p$-estimates for unimodular Fourier multipliers. Indeed, it has been shown \cite{bk:bk, stolyarov} that if $m_t(\xi)=e^{it\phi(\xi/|\xi|)}$ with $\phi\in C^\infty(\Sigma^{n-1})$, then 
\begin{equation}\label{bk}
	||e^{it\phi(\xi/|\xi|)} ||_{\mathcal M_p}\lesssim |t|^{n\left|\frac 1p-\frac12\right|}, \quad\forall t\in\mathbb R: |t|\ge 1, 
\end{equation}
thus solving Problem 15 proposed by V. Maz'ya in \cite{mazya}. In particular, the proof in \cite{bk:bk} relies on the fact that $T_{m_t}$ is essentially a Calder\'on-Zygmund operator (and hence bounded on $L^p$), and also on a precise estimate provided by A. Seeger \cite{s:s} on the weak-type boundedness of the rough operator \eqref{eq:rough}. 
\par Here, we observe  that an estimate of the form \eqref{bk} (or even not so precise) for multipliers of the form $e^{i\phi(\xi/|\xi|)}$ with a weaker condition in $\phi$ would imply the $L^p$-boundedness of the corresponding $S_\Omega$ operator associated to the multiplier $\phi(\xi')$ in the case of weighted spaces $L^p(v)$. Precisely, assume that $\Omega \in L^2(\Sigma^{n-1})$ satisfies the cancellation property. Consider the distribution 
\[ K(x)=p.v. \  \frac{ \Omega(x')}{|x|^n}, \]
and the corresponding Fourier multiplier $m_{\Omega}(x):=\widehat K (x)$. Also, let $\mathcal M(L^2(v))$ be the set of functions $m$ such that the operator $T_m$ in equation \eqref{eq:Tm} is bounded from $L^2(v)$ to $L^2(v)$.  In this context, we prove the following result: 

\begin{thm}
	\label{tomega} Fix $v$ a weight in $\mathbb  R^4$. Then, for every $\Omega\in L^2(\Sigma^3)$ such that $\int_{\Sigma^{3}}\Omega=0$ and $||\Omega||_2=1$ we have
	\begin{equation}\label{b2b2}
		S_\Omega: L^2(v) \longrightarrow L^2(v), \qquad ||S_\Omega||\lesssim 1,
	\end{equation}
	if and only if there exists $0<s<1$ so that for every $\Omega$ as before, 
	\begin{equation} \label{eq:tomega-decay}
		||e^{it \Re m_{\Omega}}||_{\mathcal M(L^2(v))} \lesssim e^{|t|^s} \quad \mbox{and}\quad ||e^{it \Im m_{\Omega}}||_{\mathcal M(L^2(v))} \lesssim e^{|t|^s}.
	\end{equation}
\end{thm}

Since a characterization of the weights for which \eqref{b2b2} holds is still unknown, we believe that the asymptotic behaviour of unimodular multipliers provides a new perspective to understand the behaviour of rough operators. 

\subsection{Singular integrals along a curve}
Consider
\begin{equation*} 
	K_{\gamma}f(x)= p.v.  \int_{\mathbb R}  {f(x-\gamma(t))} K(t) dt, \qquad x \in \mathbb R^n, 
\end{equation*}
where $ \gamma:\mathbb R \to \mathbb R^n$  is a    curve and $K$ is any kernel $K$ (if $K(t)=\frac 1t$, this is the Hilbert transform  along the curve). The $L^p$ boundedness of this operator has a long history and it is, in general, a difficult question which was first dealt with in the 1970's in a collection of papers by A. Nagel, E. Stein, N. M. Rivi\'ere and S. Wainger   \cite{stw1, stw2, sw:sw, nrw1, NRW, nrw:nrw}. 
We observe that, formally,
$$
\widehat{K_{\gamma}f}(\xi) = m_{\gamma} (\xi) \hat f(\xi),
$$
where
$$
m_{\gamma} (\xi) =  p.v. \int_{\mathbb R}  e^{i \gamma(s)\cdot \xi}K(s) ds.
$$
Hence, even if the  $L^2$ boundedness of this operator is equivalent to the boundedness of $m_\gamma$, to show the latter usually requires extensions of Van der Corput lemmata and other non-trivial techniques related to oscillatory integrals. 

Our first result in this context is the following. We emphasize here that although we initially need the kernel $K$ to be integrable, the bound of the operator will not depend on $||K||_1$. Therefore, this condition can be avoided by an approximation argument.
\begin{thm} \label{thm:curves1}
	\begin{enumerate}
		\item[i)] If, for some constant $B$, 
		\begin{equation}\label{boundedcurve}
			K_\gamma: L^p\longrightarrow L^p, \qquad ||K_\gamma|| \le  B, 
		\end{equation}
		then the functions
		\begin{equation} \label{sincos}
			\int_{\mathbb R} e^{it\cos (\gamma(s)\cdot\xi)} K(s) ds \quad\mbox{and}\quad  \int_{\mathbb R} e^{it\sin (\gamma(s)\cdot \xi)} K(s) ds
		\end{equation}
		are in $\mathcal M_p$ with norm at most $B e^{|t|}$.
		\item[ii)] Conversely, if    $K\in L^1$ and \eqref{sincos} holds with constant less than or equal $B e^{|t|^s}$ for some $0<s<1$ and some universal constant $B$, then \eqref{boundedcurve} holds with a constant independent of $||K||_1$.
	\end{enumerate}  
\end{thm}

We also prove:
\begin{thm} \label{thm:curves2}
	Let $ \gamma$ be a curve for which there exists  $1\le j\le n$ such that  $||\gamma_j||_\infty < \infty$ and 
	$$
	p.v. \int_{\mathbb R} \gamma_j(t) K(t) =\infty, 
	$$
	where $K$ is some kernel satisfying that, for every $\varepsilon>0$, 
	$$
	\int_{\varepsilon<|t|<1/\varepsilon}|K(t)| dt <\infty.
	$$
	Then, the family 
	\begin{equation} \label{uniform} K_{\gamma, \epsilon}f(x) = \int_{\varepsilon < |t| < 1/\varepsilon} {f(x-\gamma(t))} K(t) dt, \quad \varepsilon>0 \end{equation}
	is not uniformly bounded in $L^2$. 
\end{thm}

\subsection{Oscillatory integrals}

To finish, we describe an application of a different nature.
Let us consider the oscillatory integrals of convolution type given by 
$$
T_{R} f(x)=p.v.  \int_{\mathbb R} \frac{e^{iR(y)}}{y} f(x-y) dy,
$$
where $R(y)=P(y)/Q(y)$ is a rational function with real coefficients. In \cite{fw2:fw2}, it was proved that $T_R$ is bounded on $L^p$ for every $1<p<\infty$ with norm depending only on the degrees of $P$ and $Q$ but not on its coefficients. In particular, the family $(T_{tR})_{t\in \mathbb R}$ is uniformly bounded on $L^p$ . This type of uniform boundedness property  has been considered in more general cases \cite{fw:fw, fw1:fw1, p1:p1, p:p, acp:acp}.

Moreover, given two real functions $K$ and $Q$ defined in $\mathbb R^{m+n}$,  one can consider the oscillatory integral given by 
$$
T_{K, Q}f(x)= \int_{\mathbb R^n} K(x, y) e^{iQ(x, y)} f(y) dy,
$$ 
understanding the integral in the sense of principal value and $f$ to be in a test function space whenever necessary. An interesting  problem in this context is to study the asymptotic behaviour in $t\in\mathbb R$ of the boundedness in $L^p$ of the operators  $\{T_{K, tQ}\}_t$ (see, e.g.  \cite{ps:ps, b:b, p:p, g:g, gwz:gwz, ghl:ghl}). Clearly, if
\begin{equation*}\label{hkq}
	H_{KQ^n}f(x)=p.v.  \int_{\mathbb R^n} K(x, y) Q^n(x, y) f(y) dy
\end{equation*}
satisfies $||H_{KQ^n}||_p\lesssim  M^n$ for some constant $M>0$, then $||T_{K,  tQ}||_p\lesssim  e^{|t|M}$. 
We provide a weak reciprocal of this result. Although we need to assume an undesirable boundedness condition on the phase function $Q$, we shall see in Section \ref{sec:oscillatory} that this condition can be deleted on many occasions.

\begin{thm} \label{qqq}  Let $K$ be a kernel satisfying the size condition
	\begin{equation}\label{size}
		|K(x, y)| \lesssim \frac 1{|x-y|^n},
	\end{equation}
	and let $Q$ be  a bounded real function so that  there exist  constants $c, M>0$ and $0<s<1$ satisfying that, for every $z=\alpha+it\in\mathbb C$, 
	\begin{equation}\label{notanrara}
		\sup_{|x-y|<1}  \left| K(x, y) \left(e^{ z Q(x,y)}-A(z, x)\right)\right|\le  M^{|\alpha|} e^{c |t|^s}<\infty, 
	\end{equation}
	for some   analytic function $A(\cdot, x)$ with $|A(z, x)|\le M^{|\alpha|}$. Let us also assume that  there exists $c>0$ so that 
	$$
	T_{K,  tQ}: L^p\longrightarrow L^p, \qquad ||T_{K,  tQ}||\lesssim e^{c |t|^s}.
	$$
	Then, there exists $D>0$ such that 
	\begin{equation*}\label{hkq2}
		H_{KQ^n}:L^p \longrightarrow L^p
	\end{equation*}
	is bounded with constant less than or equal than $D^n$ for every $n\in \mathbb N$. 
\end{thm}

\subsection{Outline of the paper and notation}
The paper is organized as follows. In Section \ref{afo1} we shall develop our main result on Analytic families of operators proving the aforementioned extension of Stein's theorem. Section \ref{fm} is devoted to the applications to the theory of Fourier multipliers on $L^p(v)$. It contains the proofs of the results mentioned above as well as some interesting consequences, variations and improvements. Finally, Section \ref{sec:technical} deals with the study of the endpoint spaces, which is the most technical part of the paper.

Let us introduce some notation for the rest of the paper. Given two Banach spaces $A$  and $B$, we write $T: A \to B$ to indicate $T$ is \emph{bounded} from $A$ to $B$. Moreover, we write $A\approx B$ to indicate that they have
equivalent norms, $A=B$  to mean that $A\approx B$ and that the constants in the
equivalence are independent of  $\theta$, and $A\equiv B$ to indicate that $A$ and $B$ are isometrically isomorphic. 

We shall write a universal constant $C$  if $C=C(\theta)$  remains bounded when $\theta\to 0$. Please keep in mind that any universal constant may change from one occurrence to the next. As
usual, the symbol $f\lesssim g$  will indicate the existence of a positive universal constant $C$ so that $ f\le Cg$, and $f\approx g$ means that  $f\lesssim g$ and $g\lesssim f$.

\section{Analytic families of operators: the main technique}\label{afo1}
Let $\mathbb S=\{z\in \mathbb{C} : 0<\Re z<1\}$  be the unit strip and $A(\mathbb S)$ the algebra
of $\mathbb S$.  Let $\bar A=(A_0, A_1)$ be a compatible couple of Banach spaces and
$\mathcal{F}(\bar A)=\mathcal{F}(A_0,A_1)$ the space of analytic vector functions of the Calder\'on
complex interpolation method  \cite{C}; that is,   the set of all functions $f:\bar{\mathbb S}\rightarrow
A_0+A_1$ such that 

\begin{enumerate}
	\item $f$ is \emph{analytic}; that is, for every $l\in (A_0+A_1)^*$, $l(f(\cdot))\in A(\mathbb S)$,
	\item $f(z)\in A_j$ for $\text{\rm Re }z=j$ and $j=0,1$,
	\item $f(j+i\cdot)$ is $A_j$-continuous, and
	\item $\|f\|_{\mathcal{F}(\bar A) }=\max_{j=0,1}\sup_{t\in \mathbb R}\|f(j+it)\|_{A_j}<\infty$.
\end{enumerate}

The classical complex interpolation space (Calder\'on space) is defined, for $0<\theta<1$, by
$$
\bar A_\theta=\{ a=f(\theta): f\in \mathcal{F}(\bar A)\},
$$
with norm $\n a{\bar A_\theta}=\inf\{ \|f\|_ {\mathcal{F}(\bar A)}:  a=f(\theta)\}$.

In \cite{sch:sch}, the following interpolation spaces were introduced: 
$$
\bar A_{\dn}= \{a\in A_0+A_1:\, \exists f\in\mathcal{F}(\bar A), \, \fn=a\},
$$ 
with the norm
$$
\Vert a\Vert_{\dn} =\inf\{\Vert f\Vert_{\mathcal{F}(\bar A)}: \,   
\fn=a\}, 
$$
and   
$$
\bar A^{\dn}= 
\{a\in A_0+A_1\ :\ \exists   f\in\mathcal{F}(\bar A),\,  f(\t)=a,\, \fm=0,\, m=1,\ldots,n\},
$$ 
with the analog norm  $\Vert x\Vert^{\dn}$.  At this point, we should say that although all the theory could be developed for the $n$-th derivative, we shall only discuss the case $n=1$. This will be done to avoid unnecessary complications, since all our examples will be covered by this particular case.

\begin{remark} \label{rem:facts} We shall use the following classical facts repeatedly:
\begin{itemize}
	\item $\bar A^{\delta'(\theta)}\subset \bar A_\theta$ with constant $1$.
	\item $\bar A_\theta \subset \bar A_{\delta'(\theta)}$ with constant $2\theta$.
\end{itemize}
The first embedding is clear. To prove the second, 
let $\mathbb D=\{z\in\mathbb C :\vert z\vert< 1\}$ be the unit disc and consider the conformal map $\varphi_\theta:\mathbb S\longrightarrow \mathbb D
$ such that $\varphi_\theta(\theta)=0$. Explicitly, 
$$
\varphi_\t(z)=h_\alpha\bigg(\frac{e^{i\pi z}-i}{e^{i\pi z}+i}\bigg),
$$
where $\alpha=\frac{e^{i\pi \t}-i}{e^{i\pi \t}+i}$ and $h_\alpha:\mathbb D\rightarrow \mathbb D$ is a conformal
map such that $h_\alpha(\alpha)=0$. For our purposes, we need to know that
$$
\vert\varphi_\theta'(\theta)\vert = \frac \pi {2\sin \pi\t}, 
$$
and hence $\vert\varphi_\theta'(\theta)\vert\approx 1/\theta$ when $\theta\to 0$.  Now, if $a\in \bar A_\theta$ and $f\in \mathcal{F}(\bar A)$ is such that $f(\t)=a$, then the function $F(z)= \frac{\varphi_\t(z)}{\varphi'_\t(\theta)} f(z) \in \mathcal{F}(\bar A)$ and $F'(\theta)=a$. Therefore, 
$$
\|a\|_{\bar A_{\delta'(\theta)}}\le \|F\|_{\mathcal{F}(\bar A)}= \frac 1{|\varphi'_\t(\theta)|}\|f\|_{\mathcal{F}(\bar A)}\leq 2 \theta \|f\|_{\mathcal{F}(\bar A)}, 
$$
and taking the infimum over all such $f$, we obtain the result. \end{remark}

We shall repeatedly use the conformal map $\varphi_\theta$, as well as its inverse function $\phi_\theta=
\varphi_\theta^{-1}: \mathbb D\rightarrow \mathbb S$. Whenever it is clear we shall simply write $\varphi$ and
$\phi$. 

The following endpoint spaces will appear in our main theorem: 

\begin{defi} \label{defi:endpoint}
	For a given couple $(A_0, A_1)$, we consider the \emph{endpoint spaces} given by
	\begin{align*}
	\bar A_{\delta(0)} &=\bigg\{a\in A_0+A_1: \|a\|_{\delta(0)}=\overline{\lim_{\theta\to 0}} \|a\|_{\t}<\infty\bigg\}, 
	\\[2mm]
		\bar A_{\delta'(0)} &=\bigg\{a\in A_0+A_1: \|a\|_{\delta'(0)}=\overline{\lim_{\theta\to 0}}\frac{\|a\|_{\di}}{\t}<\infty\bigg\}, \\[2mm]
		\bar A^{\delta'(0)} &=\Big\{a\in A_0+A_1 : \|a\|^{\delta'(0)}=\ \overline{\lim_{\theta\to 0}}\|
		a\|^{\di}<\infty\Big\}. 
	\end{align*}
	Also, if $D$ is  a dense set 
	in the Banach space $A_0\cap A_1$ endowed with
	the norm $\|x\|_{A_0\cap A_1}=\max(\|x\|_{A_0},  \|x\|_{A_1})$, we shall denote by $D^{\delta'(0)}$ the completion of $D$ with
	respect to   $\|\cdot\|^{\delta'(0)}$, and similarly $D_{\delta'(0)}$ and $D_{\delta(0)}$.  
\end{defi}

\begin{remark} \label{rem:facts-2} As a consequence of Remark \ref{rem:facts}, we have that $D^{\delta'(0)} \subseteq D_{\delta(0)}$ with constant $1$, and $D_{\delta(0)} \subseteq D_{\delta'(0)}$ with constant $2$. 
\end{remark}

Observe that neither of the endpoint spaces $\bar{A}_{\delta(0)}$, $\bar{A}^{\delta'(0)}$, and $\bar{A}_{\delta'(0)}$ are required to be complete. Such spaces will be studied for the couples $(L^{p_0}, L^{p_1})$ --for $1\leq p_0, p_1\leq \infty$-- and  $(L^p(w_0), L^p(w_1))$ --for $1\leq p \leq\infty$ and $w_0, w_1$ two suitable weights-- in the last section of this paper. Precisely, we shall prove the following results.

\begin{claim} \label{claim:main-2} 
For every dense set $D$ in $L^{p_0}\cap L^{p_1}$, it holds that:
\begin{enumerate}
	\item The norms $\|\cdot\|_{\delta(0)}$ and $\|\cdot\|^{\delta'(0)}$ are equivalent to the $L_{p_0}$-norm on $D$. 
	\item $(L^{p_0}, L^{p_1})_{\delta(0)} \subseteq L^{p_0}$ and $(L^{p_0}, L^{p_1})_{\delta'(0)}\subseteq L^{p_0}$ up to some universal constant.
\end{enumerate}
\end{claim}

\begin{claim} \label{claim:main} Let $
	N=\{ v\  \text{\rm measurable} : v(x)\neq 0 \text{ \rm a.e. }x\}$. If $w_0$ and $w_1$ are two weights such that $w_1\in N$, then for every dense set $D$ in $L^p(w_0)\cap L^p(w_1)$, it holds that: 
\begin{enumerate}
	\item The norms $\|\cdot\|_{\delta(0)}$ and $\|\cdot\|^{\delta'(0)}$ are equivalent to the $L_{p}(w_0)$-norm on $D$.
	\item $(L^{p}(w_0), L^{p}(w_1))_{\delta(0)}\subseteq L^{p}(w_0)$ and $(L^{p}(w_0), L^{p}(w_1))_{\delta'(0)}\subseteq L^{p}(w_0)$ up to some universal constant. 
\end{enumerate}
\end{claim}

We shall now work with what we call  \emph{Calder\'on analytic families of
	operators} (we direct the interested reader to \cite{CJ, Ca1} and the references therein for more details about these families).
To simplify the computations in the applications,  we will slightly modify the classical definition.

\begin{defi}\label{afo}
	Let $\bar A$ and $\bar B$ be two compatible couples of Banach spaces, $D$ a dense subset of $A_0\cap A_1$ and 
	$\bar L=\{L_\xi\}_{\xi\in\bar{\mathbb S}}$ a family of linear operators such  that 
	$$
	L_\xi:D\longrightarrow B_0+ B_1.
	$$
	We say that $\bar L$  is a \emph{Calder\'on analytic family of operators on $D$},  and we write $\bar L\in \C(\mathbb S; D)$, if the following conditions hold:

	\begin{enumerate}
		\item For every $t\in\mathbb R$ and $j=0,1$, 
		$$
		L_{j+it}:(D,\ \Vert \cdot \Vert_{A_j})\longrightarrow (B_j,\ \Vert \cdot \Vert_{ B_j})
		$$
		is bounded with constant $M_j$.

		\item For every $a\in D$, the function $L_\xi a\in
		\mathcal F(\bar B)$.
	\end{enumerate}
	
	Moreover, if $\max(M_0, M_1)\le 1$, we say that $\bar L$ is a \emph{uniformly bounded} Calder\'on analytic family  on $D$, and we denote this fact by $\bar L\in \BC(\mathbb S,D)$.
\end{defi}

The  following  extension of Stein's theorem (see \cite{St1, CJ, CC2, Ca1})  holds. We include the proof for the sake of completeness: 

\begin{thm} \label{thm:stein}
	Let $\bar A$ and $\bar B$ be two compatible couples of Banach spaces and let  
	$\bar L\in \BC(\mathbb S; D)$. Then, the derivative operator $L'_\t=(L_\xi)'(\t)$ is bounded from $D^{\di}$ to $\bii$ with norm at most $1$.
\end{thm}
\begin{proof}
	Let 
	$$
	\mathcal{G}(\bar A)= \bigg\{
	\sum_{\text{\rm finite}}\varphi_jx_j : x_j\in A_0\cap A_1,
	\varphi_j\in A(\mathbb S)\bigg\}, 
	$$
	which is dense in $\mathcal{F}(\bar A)$ \cite{C}. Therefore, 
	the subspace of $\mathcal{F}(\bar A)$ given by
	$$
	\mathcal{G}^D(\bar A)=\bigg\{
	\sum_{\text{\rm finite}}\varphi_jx_j: x_j\in D,
	\varphi_j\in A(\mathbb S)\bigg\}, 
	$$
	is also dense in $\mathcal{F}(\bar A)$ thanks to the fact $D$ is dense in $A_0\cap A_1$. 
	Now let $d\in D$, $\varepsilon >0$ and $f\in \mathcal{F}(\bar A)$ such that $f(\t)=d$, $f'(\t)=0$ and
	$\n f{\mathcal{F(\bar A)}} \le \| d\|^{\di}+\varepsilon$.  Let us define
	$$
	H(\xi)=\frac{f(\xi)- e^{(\xi-\t)^2 }d }{\varphi_\t(\xi)^2}, 
	$$
	Then, it is clear that $H\in{\mathcal F(\bar A)}$ thanks to the hypothesis on $f$. Therefore, there exists $h\in \mathcal{G}^D(\bar A)$ such that $\| {H-h}\|_{\mathcal F(\bar A)}\le \varepsilon$.
	If we now consider $G(\xi)=h(\xi) \varphi_\t(\xi)^2+ d e^{(\xi-\t)^2 }$, we have that $G(\t)=d$, $G'(\t)=0$, $G\in \mathcal{G}^D$ and
	$\|{G-f}\|_{\mathcal F(\bar A)}\le\varepsilon$.
	Consequently, we have obtained that, for every $d\in D$, 
	$$
	\|d\|^{\di}=\inf\Big\{ \|g\|_{\mathcal F(\bar A)}: g(\t)=d, g'(\t)=0, g\in \mathcal{G}^D\Big\}.
	$$
	
	Now, let $d\in D$ and set $g\in \mathcal{G}^D$ such that $g(\t)=d$, $g'(\t)=0$ and  $\|g\|_{\mathcal{F}}\le \|d\|^{\di}+\varepsilon$.
	Let $F(\xi)= L_{\xi} g(\xi) $. Then, one can immediately see that $F\in\mathcal{F}(\bar B)$ and $\n F{\mathcal{F}(\bar B)}\le \n
	g{\mathcal{F}(\bar A)}\le \|d\|^{\di}+\varepsilon$, and since, 
	$$
	F'(\t)=L_\t' g(\t)+ L_\t g'(\t)=L_\t' d,
	$$
	we obtain that $L_\t' d\in \bar B_{\di}$ and $\n {L_\t' d}{\bar B_{\di}}\le \|d\|^{\di}+\varepsilon$. Letting $\varepsilon$
	tend to zero we obtain the result.
\end{proof}

Our first main result  in this context is the following. 

\begin{thm} \label{thm:main}
	Let $\bar A$ and $\bar B$ be two compatible couples of Banach spaces. 
	Let $\Omega=\{z\in\mathbb C : \text{\rm Re }z>0\}$ and set, for every $\xi\in\bar\Omega$, 
	$$
	T_\xi: D \longrightarrow B_0+B_1
	$$
	a linear operator such that:

	\begin{enumerate}[label=(\roman*)]
		\item For every $t\in\mathbb R$,
		$$
		T_{it}:(D,\ \Vert\cdot \Vert_{A_0})\longrightarrow B_0, \qquad ||T_{it}||\le M_0.
				$$

		\item There exists two constants $M$ and $M_1$ such that, for every $t\in\mathbb R$ and every $\alpha>0$, 
		$$
		T_{\alpha+ it}:(D,\ \Vert\cdot \Vert_{A_1})\longrightarrow B_1,  \qquad ||T_{\alpha+ it}||\le M_1 M^{\alpha}. 
		$$
				
		\item For every $\lambda>0$ and every $d\in D$, $T_{\lambda\xi}d\in \mathcal{F}(\bar B)$. 
	\end{enumerate}

	Then, 
	
	\begin{enumerate}[label=(\alph*)]
		\item The operator at the point $\xi=1$, 
		$$
		T_1: D_{\delta(0)}\longrightarrow \bar B_{\delta(0)}, \qquad ||T_1||\le M_0M.
		$$
		\item The derivative operator at the point $\xi=1$, 
		$$
		T'_1: D^{\delta'(0)}\longrightarrow   \bar B_{\delta'(0)}, \qquad ||T_1||\le M_0M(1+2|{\log M}|).
		$$
		 
	\end{enumerate}
\end{thm}

\begin{proof}
	First of all, it is clear that one can assume $M_0=1$. 
	Now, fix $0<\theta<1$ and let us consider the analytic family of operators $R_{\xi}=T_{\xi/\theta}$, where $\xi \in \mathbb S$.  Then, one can immediately see that the family $L_{\xi}=(M_1^{\t}M)^{-\xi/\theta}R_{\xi}\in
	\BC(\mathbb S; D)$. To prove (a), observe that Stein's theorem gives us that, regardless of $\theta$, 
	$$
	R_{\t}=T_1:(D, \n{\cdot}{\bar A_\theta}) \longrightarrow \bar B_\theta
	$$
	is bounded with constant $M_1^{\t} M$, and letting $\theta$ tends to zero, we clearly have that 
	$$
	T_1:  D_{\delta(0)}\longrightarrow \bar B_{\delta(0)}
	$$
	is bounded with constant $M$.
	
	\par Now, to show part (b), we apply  Theorem 2.3 to obtain that 
	$$
	L'_{\theta}: D^{\di}\longrightarrow \bar B_{\di}
	$$
	with constant $1$. Now, since 
	$$
	L'_{\theta}=\frac 1{\theta}\frac1{M_1^{\t}M}\Big( T'_1 - \log (M M_1^{\t}) T_1\Big), 
	$$ 
	letting $\theta\to 0$ we conclude that for every $d\in D$ , 
	$$
	\Big\Vert{\Big(T'_1 - \log{M}\, T_1\Big)d}\Big\Vert_{\bar B_{\delta'(0)}}\le M \Vert d\Vert^{\delta'(0)}, 
	$$
	and so we get our conclusion using part (a) and Remark \ref{rem:facts-2}. 
\end{proof}

For our applications below, we need to improve Theorem \ref{thm:main} by considering  Analytic families of operators with exponential growth at the border. To this end, we first need the following lemma.

\begin{lemma} Fix $s \in (0,1)$ and consider $\psi: \mathbb R^2 \to \mathbb R$ defined by
	\[ \psi_s(x,y) =  \frac{1}{\pi}\int_{-\infty}^{+\infty} \frac{x}{x^2+t^2}\, |y+t|^s \, dt. \]

Then, the following properties hold:
	\begin{enumerate}
		\item $\psi_s$ is harmonic. 
		\item For every $y\in \mathbb R$, $\lim_{x\to 0^+} \psi_s(x,y) = |y|^s$. 
		\item For every $x>0$, $\psi_s(x,y) \geq |y|^s$. 
	\end{enumerate}
\end{lemma}
\begin{proof} The first two properties are immediate. As for (3), we observe that 
 since $|y+t|^s  \geq |y|^s$ whenever $yt\geq0$, we have that for $y>0$,
	\[ \frac{1}{\pi}\int_{-\infty}^{\infty} \frac{x}{x^2+t^2}\,  |y+t|^s\, dt \geq \frac{1}{\pi} \int_{0}^{+\infty} \frac{x}{x^2+t^2}\, |y|^s\, dt = |y|^s. \]
	If $y<0$, the argument is similar. 
\end{proof}

\begin{thm} \label{thm:main-2}
	Let $\bar A$ and $\bar B$ be two compatible couples of Banach spaces and set, for every $\xi\in\bar\Omega$, 
	$$
	T_\xi: D \longrightarrow B_0+B_1
	$$
	a linear operator such that:

	\begin{enumerate}[label=(\roman*)]
		\item There exist constants $M_0>0$, $C>0$ and $s \in (0,1) $ such that, for every $t\in\mathbb R$,
		$$
		T_{it}:(D,\ \Vert\cdot \Vert_{A_0})\longrightarrow B_0, \qquad \|T_{it}\| \le M_0 \exp{(C|t|^s)}. 
		$$

		\item There exists two constants $M$ and $M_1$ such that, for every $t\in\mathbb R$ and every $\alpha>0$, 
		$$
		T_{\alpha+ it}:(D,\ \Vert\cdot \Vert_{A_1})\longrightarrow B_1, \qquad \|T_{\alpha+it}\| \le M_1  M^{\alpha} \exp{(C|t|^s)}. 
		$$
		 		
		\item For every $d\in D$, $T_{\xi}d\in \mathcal{F}(\bar B)$. 
	\end{enumerate}
	
	Then, the operators
	$$T_1: D_{\delta(0)}\longrightarrow \bar B_{\delta(0)}, \qquad 
	T'_1: D^{\delta'(0)}\longrightarrow   \bar B_{\delta'(0)}
	$$
	are bounded with constants depending on $M_0$, $M$, $C$ and $s$.
\end{thm}

\begin{proof} Let $\Phi_s$ be an analytic function such that $\Re \Phi_s = \psi_s$ --with the notations of the previous lemma--, and define
	\[ R_\xi = \exp{[-C\, \Phi_s(\xi)]} \ T_{\xi}, \quad \xi \in \overline{\Omega}.\]
	Let us check that the family $(R_\xi)_{\xi \in \overline\Omega}$ satisfies the hypothesis of Theorem \ref{thm:main}. Indeed, given $t>0$ we have
			\[ \|R_{it}\|_{(D, \|\cdot\|_{A_0}) \to B_0} \leq \exp{[-C \psi_s(0,t)]} \, \|T_{it}\| \leq M_0,\]
	while for any $\alpha>0$ and $t>0$, 
	\[ \|R_{\alpha+it}\|_{(D, \|\cdot\|_{A_1}) \to B_1}\leq \exp{\![-C \psi_s(\alpha, t)]} \, \|T_{\alpha +it}\|  \leq M_1 M^\alpha.\]
	Hence, part (a) of Theorem \ref{thm:main} ensures that 
	\[ R_1 = \exp{[-C \Phi_s(1)]} \ T_1\]
	is bounded from $D_{\delta(0)}$ to $\bar{B}_{\delta(0)}$ with norm less or equal than $M_0M$, and therefore, 
\[ \|T_1\|_{\bar{D}_{\delta(0)} \to \bar{B}_{\delta(0)}} \leq M_0Me^{C\psi_s(1,0)}.\]
	Additionally, applying part (b) of Theorem 2.4 we obtain that 
	\[ R'_1 = \exp{[-C\Phi_s(1)]} \left(T_1'-C\Phi_s'(1) T_1\right)\]
	is also bounded from $D^{\delta'(0)}$ to $\bar{B}_{\delta'(0)}$ with norm at most $M_0M(1+2|{\log M}|)$. This fact in tandem with Remark \ref{rem:facts-2} yields the conclusion:  
	\[ \|T'_1\|_{{D}^{\delta'(0)} \to \bar{B}_{\delta'(0)}} \leq e^{C\psi_s(1,0)} M_0M \left[1+2|{\log M}| + 2C|\Phi_s'(1)| \right]. \qedhere\]
	
	\begin{remark} The following facts will be important for applications: 
		\begin{enumerate} 
			\label{rem:bound}
			\item[\emph{(i)}] It is straightforward to see that  the conclusion of Theorem \ref{thm:main-2} also holds for $T_{\lambda+i\mu}$ and $T_{\lambda+i\mu}'$ for any fixed $w=\lambda+i\mu \in \Omega$.
	Indeed, one just need to consider $\bar T_{\xi} = T_{\lambda\xi + i\mu}$ instead of the original family $(T_\xi)_{\xi\in \Omega}$.  In other words, under the hypothesis of Theorem \ref{thm:main-2}, 	
		\[ T_\xi: D_0 \to \bar B_0, \qquad T_\xi': D^{\delta'(0)} \to \bar B_{\delta'(0)}\]
	are bounded for \emph{every} $\xi \in \Omega$. The precise bounds will be needed for the proofs of Theorems \ref{thm:curves1}, \ref{thm:curves2} and \ref{qqq} --see Sections \ref{sec:curves} and \ref{sec:oscillatory}--. We have the bounds
	\begin{align*} \|T_{\lambda+i\mu}\|_{\bar{A}_{\delta(0)} \to \bar{B}_{\delta(0)}} & \leq A(\lambda, s)\, M_0\, M\,  e^{C|\mu|^s}, \\[2mm]
	 \|T'_{\lambda+i\mu}\|_{\bar{A}^{\delta'(0)}\to \bar{B}_{\delta'(0)}} & \leq B(\lambda, s) \, M_0\, M(1+2|{\log M}|)e^{C|\mu|^s}, 
	 \end{align*}
	where $A(\lambda, s)$ and $B(\lambda, s)$ are constants depending only on $\lambda$ and $s$. 		
	
	\item[\emph{(ii)}] Observe that the norms of both $T_1$ and $T'_1$ are \emph{independent} of  $M_1$. This matter will be crucial in the proof of Theorems \ref{thm:curves1} and \ref{thm:curves2} --see Section \ref{sec:curves}.
	\end{enumerate}
	\end{remark}
	
	\begin{remark} In Theorem \ref{thm:main-2}, the dependence of the given bounds for the norms of  $T_1: D_{\delta(0)}\longrightarrow \bar B_{\delta(0)}$ and  $T'_1: D^{\delta'(0)}\longrightarrow   \bar B_{\delta'(0)}$ cannot be made independent of $s\in (0,1)$. Indeed, it can be easily checked that
	\[ \psi_s(1,0) = \frac{2}{\pi} \int_{0}^{+\infty} \frac{t^s}{1+t^2} \ dt =  \frac{1}{\pi} \int_0^{+\infty} \frac{u^{s-\frac12}}{1+u}\, du = \frac{1}{\pi} \Gamma\left(\frac{1-s}{2}\right)\Gamma\left(\frac{1+s}{2}\right). \]
	Therefore, when $s \to 1^-$, $\psi_s(1,0) \approx \frac{1}{1-s}$. The derivative $|\Phi'_s(1)|$ exhibits the same behaviour when $s \to 1^-$, since a routine argument via differentiation under the integral sign yields
	\[ |\Phi'_s(1)| = \frac{1}{\pi} \int_{-\infty}^{+\infty} \frac{u^{\frac s2+\frac12} -u^{\frac s2-\frac12} }{(1+u)^2}\, du \approx \Gamma\left(\frac{1-s}{2}\right)\Gamma\left(\frac{5+s}{2}\right).\]
	\end{remark}

\end{proof}

\par The following reformulation of Theorem \ref{thm:main-2} may be more suitable for certain applications --see the proof of Theorem \ref{sector1}--. 
In the sequel, we consider the principal argument of any complex number in $(-\pi, \pi]$. 
\begin{thm} \label{thm:main-3} Let $\bar A$ and $\bar B$ be two compatible couples of Banach spaces.	Fix $\varphi\in (0,\frac\pi2)$ and consider 
	\[ \Omega_\varphi = \{z\in \mathbb C: |{\arg(z)}|< \varphi \}.\]
	Set, for every $w\in \overline{\Omega_\varphi}$, a linear operator 
	\[ T_\xi: D \longrightarrow B_0 + B_1\]
	such that 
	\begin{enumerate}
		\item There exist constants $C>0$ and $M_0>0$ such that, for any $w\in \partial\Omega_\varphi$, 
		\[ \|T_w\|: (D, \|\cdot\|_{A_0}) \to B_0, \qquad ||T_{it}||\le M_0 e^{C|w|}. 
		 \]
				
		\item There exist constants $M_1$ and $\alpha>0$ so that, for every $w\in \Omega_\varphi$ such that $\Re(w)>\alpha$, 
		\[ \|T_w\|:  (D, \|\cdot\|_{A_1}) \to B_1, \qquad ||T_w||\le M_1 e^{C|w|}.
		\]

		\item For every $d\in D$, $T_{w}d\in \mathcal{F}(\bar B)$. 
	\end{enumerate}
		Then, the operators
	$$T_1: D_{\delta(0)}\longrightarrow \bar B_{\delta(0)}, \qquad 
	T'_1: D^{\delta'(0)}\longrightarrow   \bar B_{\delta'(0)}
	$$
	are bounded with constants depending only on $M_0$ and $C$. 
\end{thm}
	
\begin{proof}The conformal map $\rho(z) = z^{s}$, where $s = \frac{2}{\pi}\varphi<1$, sends the right half-plane $\Omega$ to $\Omega_{\varphi}$. Then, for every $z\in \overline{\Omega}$, we define $R_{z} = T_{\rho(z)}$ and check that $\bar{R} = \left(R_{z}\right)_{z\in \overline \Omega}$ satisfies the hypotheses of Theorem \ref{thm:main-2}. Indeed, using that $|\rho(z)|=|z|^s$ for any $z\in \Omega$ we obtain the following calculations: 
\begin{enumerate}
		\item For any $t\in \mathbb R$, we have:
		\[ \|R_{it}\|_{A_0 \to B_0} = \|T_{\rho(it)}\|_{A_0 \to B_0} \leq M_0 \exp{|\rho(it)|} = M_0 \exp{|t|^s}.  \]
		\item On the other hand, it can be easily checked that
		\[ \rho^{-1}(\{w\in \Omega_\varphi: \Re(w)\geq \alpha\}) \supseteq \{z\in \Omega: \Re(z) \geq \alpha^{{1}/{s}}\}, \]
		and therefore, for any $\beta > \alpha^{1/s}$ we obtain the estimation
		\begin{align*}
			\|R_{\beta + it}\|_{A_1 \to B_1} & = \|T_{\rho(\beta + it)}\| \leq M_1\exp{(|\beta+it|^s)} \leq \\
			& \leq M_1 \exp{(\beta^s + |t|^s)} \leq M_1 e^\beta \exp{|t|^s}. 
		\end{align*}
	\end{enumerate}
	Finally, it is clear that $\overline{R}$ is an analytic family, since $\rho$ is holomorphic. Hence, our previous theorem implies that $R_1: D_{\delta(0)} \to \bar{B}_{\delta(0)}$ and $R'_1: D^{\delta'(0)} \to \bar{B}_{\delta'(0)}$ with bounds depending only on $M_0$ and $C$. Since $T_1 = R_1$ and $T'_1 = \frac{1}{s}R'_1$, we conclude. 	
\end{proof}

\begin{remark}[Optimality of Theorems \ref{thm:main-2} and \ref{thm:main-3}] \label{rem:optimal} 
Before passing to applications, let us show that one cannot take $s=1$ in Theorem \ref{thm:main-2} --or, equivalently, $\varphi=\frac\pi2$ in Theorem \ref{thm:main-3}--. Therefore, the above results are optimal. Let $K(x) = \frac{1}{x}\, \chi_{[1,+\infty)}(x)$ and consider the convolution operator
\[ T(f)(x) = (f * K)(x), \]
which is well-known to be unbounded from $L_1$ to $L_1$. We define, for a given $z\in \Omega$, the linear operator
\[ T_z(f) = \left(\frac{1}{\Gamma(z+1)}  \sum_{j=1}^\infty j^{z-2} K_j \right)* f\]
where $K_j(x) = K(x) \,\chi_{[2^j, 2^{j+1}]}(x)$. Using the following estimate (see \cite{ks:ks}), which is valid for any fixed $x\in \mathbb R$ when $|y|\to +\infty$,  
\begin{equation*}
	\label{eq:gamma} 
	\frac{1}{|\Gamma(x+iy)|} \approx e^{\frac\pi2|y|}\, |y|^{\frac12-x}, 
\end{equation*}
we infer that
\begin{align*}
	\|T_{it}\|_{L_1 \to L_1} & \leq \frac{1}{|\Gamma(1+it)|} \sum_{j=1}^\infty \frac{1}{j^2} \|K_j\|_{L_1}
	\lesssim e^{\frac\pi2 |t|},   
\end{align*}
while, on the other hand, letting $b=\alpha-[\alpha]$ and $M=(2\log2)^{-1}$, 
\begin{align*}
	\|T_{\alpha + it}\|_{L_1 \to L_2} & \leq \frac{1}{|\Gamma(1+\alpha + it)|} \sum_{j=1}^\infty j^{\alpha-2} \|K_j\|_{L_2} = \\
	 & = \frac{1}{|\Gamma(1+\alpha + it)|}  \sum_{j=1}^\infty j^{\alpha-2} \cdot 2^{-(j+1)/2}   \lesssim \\
	 & \lesssim \frac{1}{|\Gamma(1+\alpha + it)|} \int_0^{+\infty} t^{\alpha-2}\cdot 2^{-t/2}\, dt \lesssim \\[2mm] 
	 & \lesssim \frac{M^\alpha}{|\Gamma(1+\alpha + it)|} \int_0^{+\infty} t^{\alpha-2}e^{-t}\, dt = M^\alpha \frac{|\Gamma(\alpha-1)|}{|\Gamma(1+\alpha+it)|} = \\[2mm]
	 &= M^\alpha \frac{(\alpha-2)(\alpha-3) \cdots |\Gamma(1+b)|}{|\alpha+it|\, |\alpha-1+it| \cdots |\Gamma(1+b +it)|} \leq \\[2mm] 
	 & \leq M^\alpha \frac{(\alpha-2)(\alpha-3) \cdots |\Gamma(1+b)|}{\alpha(\alpha-1)
	 \cdots |\Gamma(1+b +it)|}  \lesssim \frac{M^\alpha }{|\Gamma(1+b+it)|} \lesssim \\[2mm] & \lesssim M^\alpha e^{\frac\pi2 |t|}. 
\end{align*}

In light of this, observe that if our Theorem \ref{thm:main-2} could be applied to the couples $\bar{A}=(L^1, L^1)$ and $\bar{B}=(L^1, L^2)$ and the analytic family $\bar T = (T_z)_{z\in \Omega}$, then an appeal to Claim \ref{claim:main-2} would yield that $T_2(f) =\frac12 K*f$ 
is bounded from $L_1$ to $L_1$, which is false. 
\end{remark}

\section{Applications to Fourier multipliers}\label{fm}

Let us establish some notation for this section. Given $A\geq 0$ and $s\in (0,1)$, set
 \begin{align*}
 \mathcal M(A; s)&=\{ m \text{ measurable }: ||e^{itm} ||_{ \mathcal M_p}\le A e^{ |t|^{s}}, \ t\in\mathbb R\}. \\
  \mathcal M_b(A; s)&=\mathcal M(A;s)\cap L_\infty(\mathbb R^n). 
\end{align*}
Observe that, if 
 $$
 \sup_{0\le t<1}||e^{itm} ||_{ \mathcal M_p}\le A\quad\mbox{and}\quad    ||e^{inm} ||_{ \mathcal M_p}\le A    e^{ |n|^{s}}, \ n\in\mathbb Z.
 $$
 then $m\in  \mathcal M(A^2; s)$. Moreover, $m\in \mathcal M(A;s)$ does not necessarily imply that $m$ is a multiplier in $L^p$. We start the section with the proof of Theorem \ref{boundedm1}.
\begin{proof} [Proof of Theorem \ref{boundedm1}] We can clearly set $c=1$ and, hence,   $m\in  \mathcal M_b(A; s)$ for some $A>0$ and some $0<s<1$. Let us assume that $\|m\|_{\infty}>1$ (if it were not the case, we would omit the term $\|m\|_\infty$ in the definition of the following family of operators) and consider the two analytic families of operators 
\begin{equation}\label{arriba}
 T_z^{\pm} f(x)=\int_{\mathbb R^n} e^{\pm z\frac{m(\xi)}{\|m\|_\infty} }\hat f(\xi) e^{i x\xi }d\xi, \qquad z\in\Omega, 
\end{equation}
defined on the couple $\bar{A} = (L^p, L^2)$ with $D$ a test function space. It holds that
$$
\|T_{it}^\pm\|_{L^p\to L^p} \le A e^{\frac 1{\|m\|_\infty^s} {|t|^s} }\le Ae^{|t|^{s}}, \qquad \forall t\in\mathbb R,
$$
and
$$
\|T_{\alpha + it}^\pm\|_{L^2\to L^2} \le e^\alpha , \qquad \forall t\in\mathbb R, \ \alpha>0.
$$
Hence, by Theorem \ref{thm:main-2},  there exist constants $C^+_s$ and $C^-_s$ (depending only on $s$) such that 
\[ 
\|T_1^-\|_{D_{\delta(0)}\to \bar{A}_{\delta(0)}} \leq A \, C^+_s,  \quad 
\|(T_1^+)'\|_{D^{\delta'(0)} \to \bar{A}_{\delta'(0)}} \leq A \,C^-_s.\] 
Now, Claim \ref{claim:main-2} applied to the couple $(L^p, L^2)$ asserts that $T_1^-$ and $(T_1^+)'$ define multipliers on $L^p$ with
$$
\|T_1^-\|_{L^p \to L^p}=\|e^{-\frac{m}{\|m\|_\infty}} \|_{\mathcal M_p} \leq A \, C^+_s,  
$$
$$
\|(T_1^+)'\|_{L^p \to L^p}=\left |\left | \frac{m}{\|m\|_\infty}e^{\frac{m}{\|m\|_\infty}} \right |\right |_{\mathcal M_p} \leq A \,C^-_s.
$$
Therefore, their composition is also a multiplier on $L^p$ and we have that 
\begin{equation*} 
1\le  \frac{\|m\|_{\mathcal M_p}}{\|m\|_\infty}  \leq A^2 \, C^+_s \, C^-_s.  \qedhere
\end{equation*}
\end{proof}

\begin{proof} [Proof of Theorem \ref{sector1}] It is mostly analogous: define the analytic families of operators as in \eqref{arriba} and apply Theorem \ref{thm:main-3}  this time.  
\end{proof}

In another direction, the following ``interpolation" problem for multipliers is interesting: given a positive function $m\in  \mathcal M_p(\mathbb R^n)$, under which hypotheses does it hold that $m^{a}\in  \mathcal M_p(\mathbb R^n)$ for every $a>1$? Since the assertion is obvious for $a\in \mathbb N$, it is enough to consider $1<a<2$.  We provide a sufficient condition related to the behaviour of the unimodular multiplier $e^{itm}$: 
	\begin{thm} \label{thm:multiplier} 
			Let $p\geq 1$. If $m$ is a positive bounded function such that $m^{it}\in \mathcal M_p$ with $||m^{it}||_{\mathcal M_p}\lesssim e^{|t|^s}$ for some $0<s<1$, then $m^{a} (\log m)^n \in \mathcal M_p$ for every $a>0$ and every $n\in \mathbb N$. 
		\end{thm}
	
\begin{proof} Consider, for $z\in \mathbb{\overline S}$, the operator \[ T_z(f)(x) = \int m(\xi)^z\hat{f}(\xi) e^{ix\xi} d\xi. \]
Our hypothesis imply that 
\[ \|T_{it}\|_{L^p \to L^p} \lesssim e^{|t|^s}\]
while, if $\alpha>0$, the boundedness of $m$ yields
\[ \|T_{\alpha+it}\|_{L^2 \to L^2} \lesssim A\|m\|_\infty^\alpha e^{|t|^s}. \]
Hence, appealing to Theorem \ref{thm:main-2} and Claim \ref{claim:main-2}, we obtain that $T_a$ and $(T_a)'$ are bounded from $L^p$ to $L^p$ for every $a>0$; in other words, $m^a$ and $m^a \log m$ belong to $\mathcal M_p$ for every $a>0$. To obtain the result for $n>1$ we just write $m^a (\log m)^n$ as the product of $m^{a/2} (\log m)^{n-1}$ and $m^{a/2} \log m$.
\end{proof}

We can also give a couple of necessary conditions for a bounded function $m$ to be in the class $ \mathcal M(A; s)$, which we believe they are interesting in their own right.

\begin{prop} If $m\in  \mathcal M(A; s)$ for some $0<s<1$, then 
	$$
	\frac 1{\lambda +i m} \in \mathcal M_p, \qquad \forall \lambda\in\mathbb R\setminus\{0\}
	$$
	and, in general, for every $\varphi\in L^1(e^{|\cdot|^s})$ it holds that $\widehat\varphi(m)\in  \mathcal M_p$. 
\end{prop}

\begin{proof} It follows immediately by writing 
	$$
	\widehat\varphi(m(\xi))=\int_{\mathbb R} \varphi(x) e^{-i m(\xi) x}dx, 
	$$
	and applying Minkowsky integral inequality to obtain
	$$
	||\widehat\varphi(m(\cdot))||_{\mathcal M_p}\le \int_{\mathbb R} | \varphi(x)| \ ||e^{-i m(\xi) x}||_{\mathcal M_p} dx \le A  \int_{\mathbb R} |\varphi(x)| e^{|x|^s}dx \lesssim1.
	$$
	For the first part, it is enough to take $\varphi(x)=e^{-\lambda x} \chi_{(0, \infty)}$ (if  $\lambda>0$) or $\varphi(x)=-e^{-\l x}\chi_{(-\infty,0)}$ (if $\lambda<0$).
\end{proof}

\par The previous proposition, in tandem with Theorem \ref{thm:multiplier}, yield:

\begin{cor}
	If $m$ satisfies the hypotheses of Theorem \ref{thm:multiplier}, then
	\[ \frac{1}{\l\pm i\log m} \in \mathcal M_p, \qquad \forall \l \in \mathbb R\setminus\{0\}.\]
\end{cor}

\subsection{Local variations} Observe that it is not possible to avoid the boundedness of $m$ in Theorem \ref{boundedm1}, since any Fourier multiplier in $L^p$ is necessarily a bounded function. However, we can obtain some local results as follows.  

\par First, suppose that $m\in \mathcal M(A;s)$ for some $A>0$ and $0<s<1$, but $m$ is not bounded. Then,  we can partially amend this by considering
\[ \widetilde{m}(\xi) = m(\xi) \mod 2\pi, \]
which, since $m\in \mathcal M(A;s)$, it satisfies
\[ \|e^{i\widetilde mn}\|_{\mathcal M_p} \lesssim A\, e^{|n|^s}  \quad \forall n\in \mathbb Z. \]
Therefore, an application of Theorem \ref{boundedm1} shows: 
\begin{cor} If $m$ is a measurable function in the class $\mathcal M(A; s)$ such that
 \begin{equation*}
		\exists C>0 : \|e^{i\widetilde m\theta}\|_{\mathcal M_p} \leq C \qquad \forall \theta \in (0,1), 
\end{equation*}
then $\widetilde m \in \mathcal M_p$. 
\end{cor}

Our second approach stems from the well-known fact that, if $Q=Q(0, 1)$ is the unit cube in $\mathbb R^n$ and $Q_R=Q(0, R)$ then --see \cite[\S1]{ll:oo2}--  $\chi_{Q_R}\in \mathcal M_p$ and $\|\chi_{Q_R}\|_{\mathcal M_p} \leq c_p$, where $c_p$ is independent of $R$. 

\begin{thm}  \label{local} If $m \in \mathcal M(A;s)$ is locally bounded, then $m\chi_{Q_R}\in  \mathcal M_p$ and 
	$$
	\|m\chi_{Q_R}\|_\infty \le  \|m\chi_{Q_R}\|_{\mathcal M_p}\le C(A, s) c_p \,  \|m\chi_{Q_ R}\|_\infty.
	$$
\end{thm}

\begin{proof}    Since our problem is invariant under dilations, we can assume $R=1$. Let us denote $m_1 = m {\chi_{Q}}$ and assume  that $\|m_1\|_\infty \ge 1$ (otherwise, we just omit this term in the next definition of the operators).  Now,  let us consider the two analytic families of operators 
	\begin{equation*}\label{arribaa}
		T_z^{\pm} f(x)=\int_{Q} e^{\pm z \frac{m(\xi)}{\|m_1\|_\infty}} \hat f(\xi) e^{i x\xi }d\xi, \qquad z\in\Omega, 
	\end{equation*}
	on the couple $(L^p, L^2)$ as before. We have that 
	$$
	\|T_{it}^\pm\|_{L^p\to L^p} \le A\, c_p \, e^{\frac 1{\|m_1\|_\infty^s} {|t|^s} }\le A\, c_p \,  e^{|t|^{s}}, \qquad \forall t\in\mathbb R,
	$$
	and
	$$
	\|T_{\alpha + it}^\pm\|_{L^2\to L^2} \lesssim e^\alpha , \qquad \forall t\in\mathbb R.
	$$
	The result now follows from Theorem \ref{thm:main-2} as we have done in the proof of Theorem \ref{boundedm1}.
\end{proof}

Finally, we specialize to unimodular Fourier multipliers. The following result was proved in \cite[Theorem 3]{ll:oo}. 

\begin{thm} \label{edm} Let $1<p<\infty$, $p\neq 2$ and let $\varphi:\mathbb R^n\longrightarrow \mathbb R$ be a bounded measurable function satisfying that, for some constant $A>0$ and every $\lambda\in\mathbb R$, 
\begin{equation*}\label{lloo1}
\|e^{i\lambda \varphi}\|_{\mathcal M_p}\le A.\end{equation*}
 Then $\varphi$ has the following properties: 
\begin{enumerate}
	\item[(a)] There exists a closed set $E(\varphi)\subset \mathbb R^n$ of Lebesgue  measure zero such that $m$
	coincides almost everywhere with a linear function on every connected component  of the complement of $E(\varphi)$; that is, $\varphi(t)= (a_i, t)+b_i$ whenever $t\in I_i$ with $\mathbb R^n\setminus E(\varphi)=\bigcup_i I_i$.
	\item[(b)] The set of all distinct vectors $a_i$ is finite.
\end{enumerate}
\end{thm}

 As an immediate  consequence of Theorem \ref{boundedm1}, we have the following corollary  that adds two new properties to the function $\varphi$ in Theorem \ref{edm}. 
 \begin{cor} \label{edm2}
 	If $\phi$ is a function satisfying the hypothesis of Theorem \ref{edm}, then  $\varphi \in  \mathcal M_p$ and  
 $$
1\le  \frac{\|\varphi\|_{\mathcal M_p}}{\|\varphi\|_\infty}  \lesssim A^2. 
$$
 \end{cor}

\subsection{Extension to the weighted setting} 
\label{sec:weighted}
 In this section, we shall be dealing with the pair $(L^p(u), L^p(v))$ for which we have studied the endpoint spaces appearing in our theorem (see Section \ref{sec:technical}).
 
 In particular, we can (and with the same proof) extend the previous results for the class of Fourier multipliers on $L^p(v)$; that is, bounded functions $m$ so that the operator $T_m:L^p(v) \to L^p (v)$ whenever $v$ belongs to the Muckenhoupt class $A_p$. This class will be denoted by $\mathcal M_p(v)$. 
 
Let us first recall some facts on weighted theory.  Consider the \emph{Hardy-Littlewood maximal operator $M$} defined for locally integrable functions on $\mathbb R^n$ by
$$
Mf(x)=\sup_{Q\ni x}\frac{1}{|Q|}\int_Q |f(y)|dy,
$$
where the supremum is taken over all cubes $Q\subseteq\mathbb R^n$ containing $x$. It is known   \cite{m:m}  that, for every $1<p<\infty$, 
$$
M:L^p(w)\longrightarrow L^p(w) \quad\iff\quad w\in A_p,
$$
where, by definition, $w\in A_p$ if $w$ is a positive and locally integrable function such that 
$$
[w]_{A_p}:= \sup_{Q} \left(\frac 1{|Q|} \int_Q w (x) dx\right) \left(\frac 1{|Q|} \int_Q w (x) ^{1-p'}dx \right)^{p-1} <\infty. 
$$
Moreover,  if $1\le p<\infty$, 
$$
M:L^p(w)\longrightarrow L^{p, \infty}(w) \quad\iff\quad w\in A_p,
$$
where, by definition, $w\in A_1$ if there is $C>0$ such that
$$
Mw(x)\le C w(x), \quad \mbox{a.e. } x,
$$
and the infimum of  such constants $C$ is denoted by $[w]_{A_1}$. 

In this context, the fundamental result for our purpose is the quantitative version of the extrapolation theorem of Rubio de Francia \cite{r:r, dgpp:dgpp}. 

\begin{thm}
	\label{RF}
	If for some $1 \leq p_0 < \infty$ and every $w \in A_{p_0}$,
	\begin{align*} \label{ruhip}
		T: L^{p_0}(w) \rightarrow L^{p_0}(w), \quad \|T\| \le \psi([w]_{A_{p_0}}),
	\end{align*}
	where $\psi: (0,\infty) \rightarrow (0,\infty)$ is an increasing function, then for $1 < p < \infty$ and every $w \in A_p$, 
	\begin{align*}
		T: L^p(w) \rightarrow L^p(w), \quad \|T\| \le \tilde \psi([w]_{A_{p}})
	\end{align*}
	with
\[
	\tilde\psi (t) \lesssim 
		\begin{cases}
			\psi\left(C_1  t^{\frac{p_0-1}{p-1}}\right) & \text{if }\ 1 < p < p_0, \\
			\psi\left(C_2  t\right) & \text{if }\ p_0 < p < \infty,
		\end{cases} \]
	for some constants $C_1$ and $C_2$ depending on $p$ and $p_0$.
\end{thm} 

With this tool, we can now extend Theorem \ref{boundedm1} to the context of weighted Lebesgues spaces as follows.

 \begin{thm} \label{boundedmweight}Let $m$ be a bounded function, $p_0>1$  and  $v\in A_{p_0}$.   If there exists $s\in (0,1)$  such that 
\begin{equation}\label{rubio}
 \|e^{itm} \|_{ \mathcal M_{p_0}(v)}\le \Psi([v]_{A_{p_0}})e^{ |t|^{s}}, \qquad \forall t\in\mathbb R, 
\end{equation}
 for some increasing function $\Psi$, then necessarily $m\in  \mathcal M_p(v)$ for every $p>1$ and
 $$
1\le  \sup_{m\in M(A; s)} \frac{\|m\|_{\mathcal M_p(v)}}{\|m\|_\infty}= \sup_{m\in M(A, c; s)} \frac{\|m\|_{\mathcal M_p(v)}}{\|m\|_\infty}\le A([v]_{A_p}). 
 $$
 where  $A([v]_{A_p})$ is  determined depending of the relative position of $p_0$, $2$ and $p$ as follows (the constant $C$ depends on $p$ and $p_0$ and therefore it can be different in each case):
 \begin{align*}
	A([v]_{A_p}) \leq  A \left\{ 
\begin{array}{lr}
\Psi\left(C  [v]_{A_p}^{\frac{p_0-1}{p-1}}\right) & \text{if }\ p < 2 < p_0, \\[1mm]
\Psi\left(C  [v]_{A_p}\right) & \text{if }\ p_0 < 2<p, \\[1mm]
\Psi\left(C  [v]_{A_p}^{\frac{1}{p-1}}\right) & \text{if }\  p_0 < p < 2 \mbox{ or } p<p_0<2, \\[1mm]
\Psi\left(C  [v]_{A_p}^{p_0-1}\right) & \text{if }\ 2 < p_0 < p \mbox{ or } 2<p<p_0.
\end{array}\right.
	\end{align*}
\end{thm}

  \begin{proof} By the Rubio de Francia extrapolation theorem we have that \eqref{rubio} implies the existence of increasing function $\widetilde \Psi$ so that 
  $$
  \|e^{itm} \|_{ \mathcal M_2(v)}\le e^{ |t|^{s}}\widetilde\Psi([v]_{A_2}), \qquad \forall t\in\mathbb R,  \ \forall v\in A_2.
  $$
  Therefore, we can proceed as in the proof of Theorem \ref{boundedm1} 
  for the couple $(L^2(u), L^2)$ to  obtain that, for $A=\widetilde\Psi([v]_{A_2})$, 
  $$
1\le  \sup_{m\in M(A; s)} \frac{\|m\|_{\mathcal M_2(v)}}{\|m\|_\infty}= \sup_{m\in M(A, c; s)} \frac{\|m\|_{\mathcal M_2(v)}}{\|m\|_\infty}\le C_s\, A^2.
 $$
where $C_s$ is a constant depending only on $s$. Now, using Theorem \ref{RF} again we obtain the result. 
  \end{proof}

In the local case,  we can also extend Theorem \ref{local} 
since, for every $v\in A_p$,  $\chi_Q\in \mathcal M_p(v)$. 

 \begin{thm}  Let $m$  be a locally bounded function  and let $v\in A_p$.  If there exists $s\in (0,1)$  such that 
\begin{equation*}\label{rubio2}
 \|e^{itm} \|_{ \mathcal M_p(v)}\le A e^{|t|^{s}}\Psi([v]_{A_p}), \qquad \forall t\in\mathbb R
\end{equation*}
  for some increasing function $\Psi$, then necessarily $m\chi_{Q(0, R)}\in  \mathcal M_p(v)$ with 
 $$
\|m\chi_{Q(0, R)}\|_\infty \le  \|m\chi_{Q(0, R)}\|_{\mathcal M_p(v)}\le A([v]_{A_p}) \|\chi_{Q(0,1)}\|_{\mathcal M_p(v)} \|m\chi_{Q(0, R)}\|_\infty,
 $$
 with $A([v]_{A_p})$ as in Theorem \ref{boundedmweight}.
  \end{thm}

  \subsection{Homogeneous rough operators}
 In order to provide the proof of Theorem \ref{tomega}, let us give a brief refresh of definitions and  previously known results. Given $\Omega\in L^2(\Sigma^{n-1})$, we can expand  $\Omega=\sum_{j=0}^\infty Y_j$ where $Y_j$ belongs to the space of spherical harmonics  of order $j$. Moreover,    by the cancellation property of $\Omega$ we have that $Y_0=0$. Then, if we consider the distribution
$$
K(x)=p.v. \  \frac{ \Omega(x')}{|x|^n},
$$
it holds that the corresponding Fourier multiplier  $m_\Omega(x) := \widehat{K}(x)$ satisfies 
\begin{equation}\label{mimi}
 m_{\Omega}(x)= \sum_{j=1}^\infty i^{-j} \gamma_j Y_j, 
\end{equation}
where $\gamma_j\approx j^{-n/2}$  (see \cite[Theorem 4.7]{stein-weiss}).

\begin{defi}
	Given $s>0$ and $n\in \mathbb N$, the space $H^s(\Sigma_{n-1})$ is the subspace of functions $p\in L^2(\Sigma^{n+1})$ such that their spherical harmonic decomposition $p=\sum_{j=0}^\infty Y_j$ satisfies 
	\[ \sum_{j=0}^\infty (1+j)^{2s}\|Y_j\|_2^2 < \infty \]
\end{defi}

It is known \cite{dx:dx} that, for $s\in \mathbb N$, the space $H^s(\Sigma^{n-1})$ coincides with the Sobolev space $W^{s,2}(\Sigma^{n-1})$ of functions $p\in L^2(\Sigma^{n-1})$ whose all weak derivatives over the sphere up to order $s$ belong to $L^2$, endowed with the norm 
\[ \|p\|_{W^{s,2}} = \sum_{j=0}^s \|D^j_{\Sigma^{n-1}} p \|_{L^2(\Sigma^{n-1})}. \]
In fact, $\|p\|_{W^{s,2}} \approx \sum_{j=0}^\infty (1+j)^{2s}\|Y_j\|_2^2$. As a consequence, we have that 

\begin{equation*}\label{momo}
||m_{\Omega}||_{H^{n/2}(\Sigma^{n-1})} \approx ||\Omega||_2.
\end{equation*}
\par In order to avoid convergence problems in what follows,  let us define, for a fixed $l\in \mathbb N$,  
$$
\Omega_l= \sum_{j=1}^l Y_j.
$$
It is clear that $(\Omega_l)_{l=1}^\infty$ converges to $\Omega$ in $L^2(\Sigma^{n-1})$, satisfies the cancellation property and is indefinitely differentiable in $\Sigma^{n-1}$. Recall also that
$$
\Delta Y_j= -j(j+n-2) Y_j, \quad \forall j>0,
$$
and that, for a given $r>0$, the operator $\Delta^r$ is defined so that
$$
(-\Delta)^r Y_j=( j(j+n-2))^r Y_j, \quad \forall j>0.
$$

\begin{lemma} \label{omom}If $n\ge 2$, then, for every $l\in\mathbb N$, 
$$
||\Omega_l||_2 \approx ||\Delta^{n/4}_{\Sigma^{n-1}} m_{\Omega_l}||_{2}.
$$
As a consequence, if $\Omega\in C^\infty$, then 
$$
||\Omega||_2 \approx ||\Delta^{n/4}_{\Sigma^{n-1}} m_{\Omega}||_{2}.
$$
\end{lemma}

\begin{proof} 

We have that 
\begin{align*}
||\Omega_l||_2^2&=\sum_{j=1}^l ||Y_j||_2^2 = \sum_{j=1}^l j^n j^{-n} ||Y_j||_2^2 \approx \sum_{j=1}^l (j(j+n-2))^{n/2} j^{-n} ||Y_j||_2^2 
\\
&\approx \left|\left| \sum_{j=1}^l (j(j+n-2))^{n/4} i^{-j} \gamma_j Y_j \right|\right|_2^2= 
||\Delta^{n/4}_{\Sigma^{n-1}} m_{\Omega_l}||_{2}^2. \qedhere
\end{align*}
\end{proof}

We can now give the proof of Theorem \ref{tomega}.

\begin{proof}[Proof of Theorem \ref{tomega}]  Let us start by mentioning that, as a consequence of the Sobolev embedding theorem for compact manifolds \cite[Theorem 2.6]{h:h},  we have that 
\begin{equation*}\label{embedg}
H^2(\Sigma^{n-1}) \subset W^{1, p}(\Sigma^{n-1}), \qquad  \frac 12-\frac1{n-1}
<\frac1p<\frac 12.
\end{equation*}
We also have that \cite[Theorem 2.7]{h:h},   for every $q>n-1$, $W^{1, q}(\Sigma^{n-1})\subset C^0(\Sigma^{n-1})$. In particular, 
\begin{equation}\label{embed}
H^2(\Sigma^3) \subset W^{1, 4}(\Sigma^3)\subset L^\infty. 
\end{equation}

Now, given $\Omega$, let us consider the corresponding multiplier $m:=m_{\Omega}$, which is homogeneous of degree $0$. Let us expand $m$ as in \eqref{mimi} and observe that
$$
Re \ m = \sum_{k=2}^{\infty} (-1)^k  \gamma_{2k} Y_{2k},  
$$
Hence we can assume without loss of generality that $m$ is real by applying the following argument to its real and imaginary parts. Moreover, $m\in H^2(\Sigma^3)$ with  $||m||_{H^2(\Sigma^3)} \le ||\Omega||_2$,  and hence by \eqref{embed},
$$
||m||_{W^{1, 4}(\Sigma^3)} \lesssim  ||\Omega||_2,
$$

 Set $\Omega_l$ and $m_l$ as before and let us consider the homogeneous function $e^{itm_l}$, which is also  indefinitely differentiable on the sphere. Hence if we consider the multiplier
$$
m_{t,l}:=e^{itm_l} - \int_{\Sigma^{n-1}} e^{itm_l(y')} d\sigma(y'),
$$
it is known  \cite[Chapter IV, Theorem 4.7]{stein-weiss} that there exists  a function $\Omega_{t,l}\in L^2(\Sigma^{n-1})$ satisfying the cancellation property and such that  $m_{t,l}= m_{\Omega_{t,l}}$. Therefore, by hypothesis
$$
T_{\Omega_{t,l}}: L^2(v) \longrightarrow L^2(v), \qquad ||T_{\Omega_{t,l}}||\lesssim ||\Omega_{t,l}||_2.
$$
Now, since $\Omega_{t,l}$ is indefinitely differentiable, we have by  Lemma \ref{omom}
$$
  ||\Omega_{t,l}||_2^2 \lesssim ||\Delta^{n/4}_{\Sigma^{n-1}} e^{itm_l}||_{2}^2.
  $$
  Hence, since $n=4$, we have that 
  $$
  \Delta^{n/4}_{\Sigma^{3}} e^{itm_l}=\Delta_{\Sigma^{3}} e^{itm_l}= it e^{itm_l}  \Delta_{\Sigma^{3}} m_l - t^2 e^{itm_l} |\nabla_{\Sigma^{3}} m_l|^2. 
  $$
  By \eqref{mimi}, we get   that  
  $$
  ||\Delta_{\Sigma^{3}} m_l ||_2^2 \lesssim \sum_{j=1}^l  ||Y_j||_2^2\lesssim ||\Omega||_2^2.
  $$
 On the other hand,  by \eqref{embed}, 
   $$
\int_{\Sigma^{3}}|\nabla_{\Sigma^{3}} m_l|^4 \le ||m_l||_{W^{1, 4}(\Sigma^3)}^4 \lesssim ||m_l||_{H^{2}(\Sigma^3)}^4  \approx ||\Omega_l||_2^4 \lesssim ||\Omega||_2^4=1, 
$$   
and hence, for every $|t|\ge 1$, 
$$
  ||\Omega_{t,l}||_2^2 \lesssim ||\Delta^{n/4}_{\Sigma^{3}} e^{itm_l}||_{2}^2\lesssim |t|^4.
  $$
  Consequently, 
$$
\sup_{l\in \mathbb N}||e^{itm_l}||_{\mathcal M(L^2(v))} \lesssim (1+|t|^2) \lesssim e^{|t|^s},
$$
and the result is obtained by letting $l\to\infty$.

\par The converse implication follows from an argument similar to the proof of Theorem \ref{boundedm1}; let us sketch it. Define the corresponding analytic families for $\Re m_\Omega$ and $\Im m_\Omega$ as in \eqref{arriba}, and use the hypotheses to verify that the conditions of Theorem \ref{thm:main-2} are satisfied. Then, apply Claim \ref{claim:main} to conclude that $m_\Omega \in \mathcal M_2$, which is equivalent to the boundedness of $S_{\Omega}: L^2(v) \to L^2(v)$ with $\|S_{\Omega}\| \lesssim 1$. 
 \end{proof}

\begin{remark}

\begin{enumerate}
	\item[\emph{(i)}] Observe that we have proved that, for every homogeneous function of order $0$, $m_\Omega$,  with $\Omega$ satisfying the hypothesis of the previous theorem, we have
	\begin{equation}\label{multasy}
		\left\|e^{itm_{\Omega}}-  \int_{\Sigma^{n-1}} e^{itm_\Omega(y')} d\sigma(y')
		\right\|_{\mathcal M(L^2(v))} \lesssim \max( |t|, |t|^2), \qquad \forall t\in\mathbb R.
	\end{equation}
	\noindent In fact, we have a self-improving condition, since \eqref{b2b2} implies a much better decay in $t$ than the one given in \eqref{eq:tomega-decay}. Precisely: 
\begin{equation} \label{eq:vrf}
	||e^{it Re\, m_{\Omega}}||_{\mathcal M(L^2(v))} \lesssim 1+|t|^2 \quad \mbox{and}\quad ||e^{it Im\, m_{\Omega}}||_{\mathcal M(L^2(v))} \lesssim 1+|t|^2.
\end{equation}

\item[\emph{(ii)}] By \eqref{tomega1}, we have that \eqref{multasy} holds for every $v\in A_1$. Therefore, if one is able to prove \eqref{eq:vrf} for every $v\in A_2$, then the Rubio de Francia extrapolation theorem allows to pass to the weighted boundedness on $L^p(v)$ for every $v\in A_p$.
\end{enumerate}
\end{remark}

\subsection{Singular operators on curves} 
\label{sec:curves}
Let us provide now the proofs of Theorems \ref{thm:curves1} and \ref{thm:curves2}. Here we demonstrate the full strength of the derivative operator technique by making use of Remark \ref{rem:bound}.

\begin{proof}[Proof of Theorem \ref{thm:curves1}] 
We first deal with \emph{(i)}. If \eqref{boundedcurve} holds, then the function
$$
m_{\gamma} (\xi) =  \int_{\mathbb R} e^{i \gamma(u)\cdot \xi}K(u)\, du
$$
clearly satisfies that 
$$
\sup_{t\in\mathbb R} ||m_{t\gamma}||_{\mathcal M_p}\le B. 
$$
From here, it follows that, for every $n\in\mathbb Z$, the functions
\[ a_n(\xi) = \int_{\mathbb R}\cos (n(\gamma(u)\cdot \xi) ) K(u) du\]
are in $\mathcal M_p$ \emph{uniformly} in $n$ with norm less than or equal to  $B$. Now, using the identity
\begin{equation} \label{eq:cos}
	(\cos x)^k = \frac{1}{2^k} \sum_{j=0}^n \binom{k}{j}\cos[(k-2j)x]
\end{equation}
the function
\[\int_{\mathbb R} (\cos (\gamma(u) \cdot \xi))^k K(u) du, \quad k\in\mathbb N,  
\]
also belongs to $\mathcal M_p$ with norm at most $B$. 
Hence, by a Taylor expansion, the first function in \eqref{sincos} is a Fourier multiplier on $L^p$ with norm less than or equal $B e^{|t|}$. An analogous reasoning provides boundedness in $\mathcal M_p$ of the second function in \eqref{sincos}. 

\par To show \emph{(ii)}, we shall again proceed with the cosine factor. The work for the sine factor is analogous. Let us define the AFO
$$
T_z  f(x)=\left[ \left( \int_{\mathbb R} e^{z\cos (\gamma(u)\cdot\xi)} K(u)\, du\right) \hat f(\xi)\right]^{\vee}(x).
$$
Then, by assumption 
$$
T_{it}:L^p\to L^p, \qquad ||T_{it}||\le B e^{|t|^s},
$$
and if $z=\alpha +it$ we have that 
$$
\left| \int_{\mathbb R} e^{(\alpha +it) \cos (\gamma(u)\cdot\xi)} K(u) \, du\right| \le e^{\alpha} ||K||_1,
$$
and hence
$$
T_{\alpha+it}:L^2\to L^2, \qquad ||T_{\alpha+it}||\lesssim e^{\alpha}\|K\|_1 \leq e^{\alpha} e^{|t|^s} ||K||_1. 
$$
Applying now our Theorem \ref{thm:main-2},  we obtain that $T'_{1+it}:L^p\to L^p$ with a constant independent of $\|K\|_1$. Indeed, Remark \ref{rem:bound} informs us that 
\[|| T'_{1+it}||_{L^p \to L^p} \lesssim B e^{|t|^s}.\] 
Since the multiplier associated to this operator is given by 
$$
\int_{\mathbb R} e^{(1+it)\cos (\gamma(u)\cdot\xi)} \cos (\gamma(u)\cdot\xi) K(u)\, du, 
$$
we must now eliminate the factor $e^{(1+it)\cos (\gamma(u)\cdot\xi)}$ to conclude. For this purpose, we consider a new AFO given by 
\[ S_z f(x) = \left[ \int_\mathbb R \cos(\gamma(u)\cdot \xi)e^{(1-z)\cos(\gamma(u)\cdot\xi)} K(u)\, du\right]^\vee(x) \]
The above calculations imply that 
\[ S_{it}: L^p \to L^p, \qquad \|S_{it}\|\leq B_s e^{|t|^s}, \]
while 
\[ S_{\alpha+it}: L^2 \to L^2, \qquad \|S_{\alpha+it}\| \lesssim \|K\|_{L^1} e^{\alpha}. \]
Hence $S_{1}$ is bounded from $L^p$ to $L^p$. Therefore, its associated multiplier
\[ \int_\mathbb R \cos(\gamma(u)\cdot \xi) K(u)\, du\]
belongs to $\mathcal M_p$. The same procedured can be applied to the sine function, and so Theorem \ref{thm:curves1} follows inmediately. 
\end{proof}

\begin{proof}[Proof of Theorem \ref{thm:curves2}]
	Let us take the contrapositive route and assume that the family \eqref{uniform} is uniformly bounded. Set, for each fixed $\varepsilon>0$,  the AFO
	$$
	T_z (f)=  \left( \int_{\varepsilon<|u|<1/\varepsilon}  e^{z \gamma_j(u)} e^{i \overline \gamma(u)\cdot \bar \xi}K(u)\, du\right) \hat f, 
	$$
	where we use the standard notation $\bar y=(y_1, \cdots, y_{j-1}, y_{j+1}, \cdots, y_n)$. 
	By assumption, the family of multipliers
	\[m_{\gamma, \varepsilon} (\xi) =  \int_{\varepsilon<|u|<1/\varepsilon} e^{i \gamma(u)\cdot \xi}K(u)\, du \]
	is uniformly bounded by some constant $A$, say.
	Therefore, writing $(t, \bar \xi)= (\xi_1, ..., \xi_{j-1}, t, \xi_{j+1}, ..., \xi_n)$ and using the fact that $m_{\gamma, \varepsilon}$ is continuous, we have 
	\[ ||T_{it}||_{L^2\to L^2} = \|m_{\gamma,\varepsilon}(t, \cdot)\|_{\infty} \leq A, \qquad \forall t\in\mathbb R.\]
	On the other hand, 
	\begin{align*}
	\|T_{\alpha + it}\|_{L^2 \to L^2} & \leq \left|\int_{\varepsilon<|u|<1/\varepsilon} e^{\alpha \gamma_j(u)}\, e^{i\bar\gamma(u)\cdot (t, \bar\xi)} K(u)\, du \right| \leq  \\[2mm]
	& \leq e^{\alpha \|\gamma_j\|_\infty} \left(\int_{\varepsilon<|u|<1/\varepsilon} |K(u)|\, du \right).
	\end{align*}
	Consequently, by Theorem \ref{thm:main-2} and Remark \ref{rem:bound}
	$$
	T'_{1+it}f= \left( \int_{\varepsilon<|u|<1/\varepsilon}  \gamma_j(u) e^{(1+it) \gamma_j(u)} e^{i \overline \gamma(u)\cdot \bar \xi} K(u)\, du\right) \hat f,
	$$
	is $L_2$-bounded with a constant independent of $\varepsilon$. We now proceed as in the previous theorem to eliminate the exponential factor $e^{(1+it) \gamma_j(u)}$. This way we obtain that, taking $\overline\xi=0$, 
	$$
	\left|\int_{\varepsilon<|u|<1/\varepsilon}    \gamma_j(u) K(u) \, du\right| \lesssim A e^{||\gamma_j||_\infty},
	$$
	and letting $\varepsilon\to 0$ we obtain the result. 
\end{proof}

To conclude, observe that our Theorem \ref{thm:curves2} provides a necessary condition for the boundedness of the corresponding maximal operator (cf.  \cite{KLO22}): 
\begin{cor} Under the hypothesis of Theorem \ref{thm:curves2}, the maximal operator 
	\[ M_\gamma f(x) = \sup_{\varepsilon>0} \left|\int_{\varepsilon<|u|<1/\varepsilon} f(x-\gamma(u))\, K(u)\, du\right| \]
is unbounded from $L^2$ into $L^2$. 
\end{cor}

\subsection{Oscillatory integrals}
\label{sec:oscillatory} 

We provide the proof of Theorem \ref{qqq}. To alleviate notation, we consider the following  definition: 

\begin{defi} Given $s\in (0,1)$, we say that $\{T_{K,  tQ}\}_t$ satisfies the \emph{$s$-property on $L^p$} if there exists $c>0$ so that 
	$$
	T_{K,  tQ}: L^p\longrightarrow L^p, \qquad ||T_{K,  tQ}||\lesssim e^{c |t|^s}.
	$$
\end{defi}

\begin{proof}[Proof of Theorem \ref{qqq}.] For simplicity, we shall assume $c=1$ and by linearity, we can assume that $||Q||_\infty=1$. Let us define the AFO 
	$$
	Q_{z}f(x)= T_{z} f(x)-A(z, x) T_0f(x),
	$$
	where
	$$
	T_zf(x)=p.v.  \int_{\mathbb R^n}  K(x, y) {e^{z Q(x, y)}} f(y) dy, \qquad \Re z>0. 
	$$
	By hypothesis, there exists $M_0>0$ such that 
	$$
	T_{it}: L^p \longrightarrow L^p, \qquad ||T_{it}||\le M_0 e^{ |t|^s}, \ \forall t\in\mathbb R, 
	$$
	and hence 
	$$
	Q_{it}: L^p \longrightarrow L^p, \qquad ||Q_{it}||\le 2M_0 e^{ |t|^s}, \ \forall t\in\mathbb R, 
	$$
	On the other hand,  let us write 
	\begin{align*}
		Q_{z}f(x)&=p.v.  \int_{|x-y|>1} K(x, y) {{e^{z Q(x, y)}}} f(y) dy 
		\\
		&+   p.v. \int_{|x-y|<1}  K(x, y)(e^{z Q(x, y)}-A(z, x)) f(y) dy 
		\\
		&-A(z,x)  \int_{|x-y|>1} K(x, y) f(y) dy.
	\end{align*}
	
	It  is immediate to see that, for every $\alpha>1$, 	
	$$
	Q_{\alpha+it}: L^1\longrightarrow L^\infty, \qquad ||Q_{\alpha+it}||\le (e^{\alpha}+ 2M^\alpha)e^{|t|^s}\le (e+2M)^\alpha e^{|t|^s}
	$$
	and so an appeal to Theorem \ref{thm:main-2} and Claim \ref{claim:main-2} yields 
	$$
	Q_{1+it}': L^p\longrightarrow L^p, \qquad ||Q_{1+it}'||\le B_s M_0 (e+2M)^2  e^{|t|^s},
	$$
	with $B_s$ depending only on $s$.   Set $k= \sup_{t, x} |A'(\pm 1+it, x)|$. Then, the operator
	$$
	p.v.  \int_{\mathbb R^n} Q(x,y) e^{ Q(x, y)}K(x, y)e^{it Q(x, y)} f(y) dy= Q_{1+it}'f(x)+ A'(1+it, x) T_0f(x)
	$$
	is bounded from $L^p$ to $L^p$ with norm less than or equal to
	$
	B_s M_0 (e+2M)^2(1+k)e^{|t|^s}.
	$
	
	Now, in the case $n=1$, we define  the AFO
	\begin{align*}
		Q_{z, 1}f(x)&=p.v.  \int_{|x-y|>1} Q(x,y) e^{ Q(x, y)}  K(x, y) {{e^{-z Q(x, y)}}} f(y) dy 
		\\
		&+   p.v. \int_{|x-y|<1}  Q(x,y)  e^{ Q(x, y)} K(x, y)(e^{-z Q(x, y)}-A(-z, x)) f(y) dy 
		\\
		&-A(-z, x) \int_{|x-y|>1} Q(x,y)e^{Q(x,y)}K(x,y) f(y) dy.
	\end{align*}
	Then, it holds that 
	$$
	Q_{it, 1}:L^p \longrightarrow L^p, \qquad ||Q_{it, 1}||\le B_s M_0 (e+2M)^2(1+k)e^{|t|^s}.
	$$
	while
	$$
	Q_{\alpha+it, 1}: L^1\longrightarrow L^\infty, \qquad ||Q_{\alpha+it}||\le e(e^{\alpha}+ M^\alpha)e^{|t|^s}\le e (e+2M)^\alpha e^{|t|^s},
	$$
	 
	Therefore, we obtain that $ Q_{1, 1}:L^p\longrightarrow L^p$ with $||Q_{1,1}||\le B_s^2 M_0 (e+2M)^3(1+2k)$. Since 
	$$
	Q_{1, 1}f(x)=p.v.  \int_{\mathbb R^n} K(x, y){Q(x, y)}  f(y) dy- A(-1+it, x)T_0f(x),
	$$
	the result follows with constant less than or equal to $B_s^2 M_0 (e+2M)^4(1+2k)$.
	
	In the case $n=2$, we take
	$$
	Q_{1, 1}'f(x)=-p.v.  \int_{\mathbb R^n} K(x, y)Q(x, y)^2  f(y) dy- A'(-1+it, x)T_0f(x), 
	$$
	with $||Q'_{1,1}||\le B_s^2 M_0 (e+2M)^4(1+2k)$ and hence 
	$$
	||H_{K Q^2}||\le B_s^2 M_0 (e+2M)^4(1+3k).
	$$

	By induction, we get that 
	$$
	||H_{KQ^n}||\le B_s^n M_0 (e+2M)^{4n}  (1+(n+1)k), 
	$$
	and the result follows. 
\end{proof}

\begin{remark} 
	\begin{enumerate}
		\item[(a)] Observe that Theorem \ref{qqq} does not hold if $Q$ is unbounded, as the simple example  $K(y)=\frac 1y$ and $Q(y)=y$ shows.  However, as we shall see in the next corollary, we can, essentially, reduce to this case by considering sinusoidal phases. 
	
\item[(b)] Let us analyze condition \eqref{notanrara} in some particular cases. Assume that $K(x, y)=\frac 1{x-y}$ and $Q(x, y)=\varphi(x-y)$ with $\varphi$ bounded and derivable at $0$. Then, if we take $A(z, x-y)=e^{ z \varphi(0)}$, we have that
	\begin{align*}
		\sup_{|x-y|<1}  &\left| K(x, y) \left(e^{ z Q(x,y)}-A(z, x)\right)\right|=\sup_{|y|<1}  \left| \frac 1y \left(e^{\pm z \varphi(y)}-e^{ z\varphi(0)}\right)\right|
		\\[2mm]
		&=\sup_{|y|<1} \left|\frac{\varphi(y)-\varphi(0)}y e^{\pm z \varphi(0)} \frac{e^{\pm z (\varphi(y)-\varphi(0))}-1}{\varphi(y)-\varphi(0)}\right|, 
	\end{align*}
	and  \eqref{notanrara}  follows. In $n$ dimensions, if $\varphi$ is radial and its radial part is bounded and derivable at $0$, the result also easily follows. The general idea is that $A(z,x)$ is of the form $e^{z T_n(Q(x, \cdot)(0)}$ where $T_n$ denotes a suitable Taylor polynomial of order $n$. 
	\end{enumerate} 
\end{remark}

The following result allows to reduce to the case of bounded phases on certain occasions. 

\begin{cor}  \label{cor:bdd-phase} Let  $K$ be a kernel satisfying \eqref{size} and such that there exists $M>0$ satisfying 
	\begin{equation*}\label{notanrara:sin}
		\sup_{|x-y|<1}  \left| K(x, y) \left(e^{ z \sin Q(x,y)}-A_1(z, x)\right)\right|\le  M^{|\alpha|} e^{|t|^s}<\infty, \qquad z=\alpha+it.
	\end{equation*}
	and
	\begin{equation*}\label{notanrara:cos}
		\sup_{|x-y|<1}  \left| K(x, y) \left(e^{ z \cos Q(x,y)}-A_2(z, x)\right)\right|\le  M^{|\alpha|} e^{|t|^s}<\infty, \qquad z=\alpha+it.
	\end{equation*}
	for some   analytic functions $A_j(\cdot, x)$ with $|A_j(z, x)|\le M^{|\alpha|}  $.
	  
	\begin{enumerate}
		\item[(i)] If $\{T_{K, t \sin Q}\}_t$ and $\{T_{K, t \cos Q}\}_t$ satisfy the $s$-property for some $0<s<1$, then
	$$
	T_{K, Q}: L^p\longrightarrow L^p.
	$$
	
		\item[(ii)] If $\{T_{K, tQ}\}_t$ satisfies  the $s$-property with $s=1$,  then so do $\{T_{K, t \sin Q}\}_t$ and $\{T_{K, t \cos Q}\}_t$.
	
	\item[(iii)]  If $\{T_{K, t \sin rQ}\}_t$ and $\{T_{K, t \cos rQ}\}_t$ satisfy the $s$-property for some $0<s<1$ uniformly in $r\in\mathbb R$, then 
	$\{T_{K, t Q}\}_t$ satisfies the uniform bound in $t\in\mathbb R$. Consequently, $\{T_{K, t \sin rQ}\}_t$ and $\{T_{K, t \cos rQ}\}_t$ satisfy the $s$-property for $s=1$ uniformly in $r\in\mathbb R$.
	\end{enumerate}
\end{cor}

\begin{proof} Item \emph{(i)} follows easily from Theorem \ref{qqq}, since it asserts that $H_{K\sin Q}$ and $H_{K\cos Q}$ are bounded in $L^p$ and
	$$
	T_{K, Q}=H_{K\cos Q} + i H_{K\sin Q}.
	$$
	
	To prove \emph{(ii)},  we appeal once again to the formula \eqref{eq:cos}. Hence
	$$
	H_{K (\cos \, Q)^n} = Re\left( \frac 1{2^n}\sum_{k=0}^n  \binom{n}{k}\ T_{K, (n-2k)Q}\right). 
	$$
	Now, by hypothesis, there exists $A>0$ such that $||T_{K, (n-2k)R}||\le Ae^{|n-2k|}$, which in turn implies
	\begin{align*}
		\|H_{K (\cos \, Q)^n}\|& \le \frac 1{2^n}\sum_{k=0}^n  \binom{n}{k}\ A^{|n-2k|} =
		\\[2mm]
		&=  \frac 1{2^n}\sum_{k=0}^{[n/2]}  \binom{n}{k}\ A^{n-2k}+  \frac 1{2^n}\sum_{k=[n/2]+1}^n  \binom{n}{k}\ A^{2k-n} \leq D^n,
	\end{align*}
	where $D=\frac{A^2+1}{A}$. Therefore,  since
	$$
	T_{K, t \cos Q}= \sum_{n=0}^\infty \frac {(it)^n}{n!}  H_{K (t\cos Q)^n},
	$$
	we obtain that $||T_{K, t \cos Q}\| \le e^{D|t|}$ and thus  $\{T_{K, t \cos Q}\}_t$ satisfies the $1$-property. A similar reasoning works for $\{T_{K, t \sin Q}\}_t$. Finally, item \emph{(iii)} follows by arguing as in \emph{(i)} and \emph{(ii)} keeping track of the constants. 
\end{proof}

\begin{remark}
	If $R=P/Q$ is a rational function and, for some $0<s<1$, 
	$$
	||T_{K,  t\sin R}||\le C e^{|t|^s} \quad\mbox{and}\quad  ||T_{K,  t\cos R}||\le C e^{|t|^s}
	$$
	with $C$ only depending on the grades of $P$ and $Q$ and not on the coefficients, then $||T_{K, R}||$ is also independent on the coefficients of $P$ and $Q$ and thus $||T_{K, tR}||$ satisfies the $1$-property uniformly on $t\in \mathbb R$. In such a case, an application of Corollary \ref{cor:bdd-phase} asserts that $\{T_{K, t \cos R}\}_t$ and $\{T_{K, t \sin R}\}_t$ satisfy the $1$-property. 
	
\end{remark}

\section{The endpoint spaces: the technical part}

\label{sec:technical}
The purpose of this section is to give a proof of Claims \ref{claim:main-2} and \ref{claim:main}.

\subsection{The couple $(L^{p_0}, L^{p_1})$}
 
Let us start recalling a few items which will be useful for us here. As usual, given $\theta \in (0,1)$, define $p(\theta)$ by the equality 
\[\frac{1}{p(\theta)} = \frac{1-\theta}{p_0} + \frac{\theta}{p_1}\]
with the obvious interpretation if $p_1=\infty$. It is well-known that $(L^{p_0}, L^{p_1})_{\theta} = L^{p(\theta)}$. Hence, a standard convergence argument yields:
\begin{prop}
	Let $D$ be a dense subset of $L^{p_0} \cap L^{p_1}$. Then the norms $\|\cdot\|_{\delta(0)}$ and $\|\cdot\|_{p_0}$ agree on $D$. In particular, the inclusion ${(L^{p_0}, L^{p_1})}_{\delta(0)} \hookrightarrow L^{p_0}$ is bounded. \end{prop}

\par Concerning the Schechter interpolation spaces, we have that, for a given $\theta \in (0,1)$: 
\begin{itemize}
	\item \cite[Theorem 2.4]{C03} A function $f$ belongs to the space $(L^{p_0}, L^{p_1})^{\delta'(\theta)}$ if and only if $f(1+|{\log f}|) \in L^{p(\theta)}$, and 
	\[ \|f\|^{\delta'(\theta)} \approx \begin{cases} \|f\|_{p(\theta)} + \dfrac{\theta|{p_1-p_0}|}{(1-\theta)p_1 + \theta p_0} \left\| f \log\dfrac{|f|}{\|f\|_{p(\theta)}}\right\|_{p(\theta)} & \text{ if } p_1 < \infty,  \\[6mm]
		\|f\|_{p(\theta)} + \dfrac{\theta}{1-\theta}\left\| f \log\dfrac{|f|}{\|f\|_{p(\theta)}}\right\|_{p(\theta)} & \text{ if } p_1 = \infty. \end{cases}
	 \]
	\item \cite[Theorem 2.7]{C03} Consider the space 
	\[ M_{\psi_\theta} = \left\{ f: f =\frac{f_0}{\theta} + f_1 \log \frac{|f_1|}{\|f_1\|_{p(\theta)}}; \ f_0, f_1 \in L^{p(\theta)}\right\}, \]
	endowed with the norm 
	\[ \|f\|_{M_{\psi_\theta}} = \inf\{ \|f_0\|_\theta + \|f_1\|_\theta\}. \]
	Then, $(L^{p_0}, L^{p_1})_{\delta'(\theta)} \approx M_{\psi_\theta}$, that is to say,  $(L^{p_0}, L^{p_1})_{\delta'(\theta)}$ and $M_{\psi_\theta}$ are isomorphic with constants independent of $\theta$. 
\end{itemize}

\par

\begin{thm} Let $1\leq p_0, p_1 \leq +\infty$, and let $D$ be a dense subset of $L^{p_0} \cap L^{p_1}$. Then:
	\begin{enumerate}
		\item[i)] The norms $\|\cdot\|^{\delta'(0)}$ and $\|\cdot\|_{p_0}$ are equivalent on $D$.
		\item[ii)] The inclusion ${(L^{p_0}, L^{p_1})}_{\delta'(0)} \hookrightarrow L^{p_0}$ is bounded. 
	\end{enumerate}
	\begin{proof}
		Item (i) was shown in \cite[Theorem 2.5]{C03}, and it actually follows directly from the above considerations. 
		
	To show (ii), let $f\in {(L^{p_0}, L^{p_1})}_{\delta'(0)}$, choose $(\theta_k)_{k=1}^\infty$ any decreasing sequence in $(0,1)$ converging to $0$ and let $p_k = p(\theta_k)$, which clearly converges to $p_0$. We can assume without loss of generality that $f\neq 0$ and that $f\in M_{\psi_{\theta_k}}$ for all $k\in \mathbb N$ --eliminating finite terms from the sequence $(\theta_k)_{k=1}^\infty$ if necessary--. Now, given $\varepsilon>0$, for every $k\in \mathbb N$ there is a decomposition 
		\begin{equation} \label{eq:decomposition} f = \frac{f_{0,k}}{\theta_k} + f_{1,k} \log\frac{|f_{1,k}|}{\|f_{1,k}\|_{p_k}}
		\end{equation}
		such that 
		\[ \|f\|_{M_{\psi_{\theta_k}}} \leq \|f_{0,k}\|_{p_k} + \|f_{1,k}\|_{p_k} \leq \|f\|_{M_{\psi_{\theta_k}}} + \varepsilon \theta_k\]
		Hence
		\[ \limsup_k \frac{\|f_{0,k}\|_{p_k} + \|f_{1,k}\|_{p_k}}{\theta_k} \leq \limsup_k \frac{\|f\|_{M_{\psi_{\theta_k}}}}{\theta_k} +\varepsilon < \infty.\]
		In particular, 
		\begin{equation} \label{eq:bounded} \limsup_k \frac{\|f_{1,k}\|_{p_k}}{\theta_k} < \infty, \end{equation} 
		hence $\lim_k \|f_{1,k}\|_{p_k} = 0$, and an application of Fatou's lemma yields
		\begin{equation} \label{eq:0-ae} \liminf_k |f_{1,k}| = 0 \quad a.e. \end{equation}

		Now, observing the decomposition  
		\begin{align*} f_{1,k} \log\frac{|f_{1,k}|}{\|f_{1,k}\|_{p_k}} = 
		 f_{1,k} \log{|f_{1,k}|} + f_{1,k} \log \frac{\theta_k}{\|f_{1,k}\|_{p_k}} + \frac{f_{1,k}}{\theta_k} \theta_k \log \frac{1}{\theta_k}, \end{align*}
		it is easy to deduce by means of \eqref{eq:bounded} and \eqref{eq:0-ae} that
		 \[ \liminf_k \left(f_{1,k} \log\frac{|f_{1,k}|}{\|f_{1,k}\|_{p_k}}\right)= 0 \quad a.e.\]
		 But this fact, in tandem with \eqref{eq:decomposition}, implies that
		 \[ \liminf_k \frac{f_{0,k}}{\theta_k} = f \quad a.e., \]
		 and therefore 
		 \begin{align*} \int |f(x)|^{p_0} \ dx & \leq \int \liminf_k \left|\frac{f_{0,k}(x)}{\theta_k}\right|^{p_k}\, dx \leq \liminf_k \int \left|\frac{f_{0,k}(x)}{\theta_k}\right|^{p_k}\, dx \leq 
		 \\[2mm]
		 & \leq \limsup_k \left(\frac{\|f_{0,k}\|_{p_k}}{\theta_k}\right)^{p_k} \leq \limsup_k \left(\frac{\|f_{0,k}\|_{p_k} + \|f_{1,k}\|_{p_k}}{\theta_k}\right)^{p_k} =  \\[2mm]
		  &  =  \left[\limsup_k \left(\frac{\|f_{0,k}\|_{p_k} + \|f_{1,k}\|_{p_k}}{\theta_k}\right)\right]^{p_0} \leq \\[2mm]
		  & \leq \left[ \limsup_k \frac{\|f\|_{M_{\psi_{\theta_k}}}}{\theta_k} + \varepsilon\right]^{p_0} < +\infty. \end{align*}
		 In other words, $f\in L^{p_0}$, and letting $\varepsilon\to 0$, we obtain that
		 \[ \|f\|_{p_0} \leq \limsup_{\theta \to 0} \frac{\|f\|_{M_{\psi_\theta}}}{\theta}  \approx \|f\|_{\delta'(0)}. \qedhere \]
	\end{proof}
\end{thm}

\subsection{The couple $(\bC, \lambda\bC)$}
To obtain the desired information regarding the endpoint spaces for the couple $(L^p(u), L^p(v))$ we need to start with the very simple one.  It is well-known that, for the classical Calderón interpolation spaces, 
\[ (\bC, \lambda\bC)_\theta \equiv \lambda^\theta\bC.\]
Our purpose is to analyze the Schechter interpolation spaces $\bar{A}_{\delta'(\theta)}$ and $\bar{A}^{\delta'(\theta)}$ for this couple. 
\par To do so, let us recall a few standard facts.  Given $f \in \mathcal F$ and $0<\theta<1$, it is known \cite[Lemma 4.3.2]{bl:bl} that 
\[ \|f(\theta)\|_\theta \leq \exp\left(\, \sum_{j=0,1} \int_{-\infty}^{+\infty} \log \|f(j+it)\|_j \, \mu_j(\theta, t)\, dt\right), \]
where $\mu_j$, $j=0,1$, are the Poisson kernels on the strip. From the above inequality plus Jensen's inequality (twice) it follows that for every $p\geq 1$,
\[ \|f(\theta)\|_\theta \leq \left(\frac{1}{1-\theta}\right)^\frac{1-\theta}{p} \left(\frac{1}{\theta}\right)^\frac{\theta}{p} \left[\, \sum_{j=0,1} \left(\int_{-\infty}^{+\infty} \|f(j+it)\|^p_j \, \mu_j(\theta, t)\, dt\right)\right]^\frac1p. \]
Therefore,  we define $\mathcal F_{\theta,p} = \mathcal F$ endowed with the norm 
\[ \|f\|_{\mathcal F_{\theta,p}} = \left[\, \sum_{j=0,1} \left( \int_{-\infty}^{+\infty} \|f(j+it)\|_j^p \, \mu_j(\theta, t)\, dt\right) \right]^\frac1p.\]
Observe that, for $1\leq p\leq q < \infty$,  $\mathcal F \subseteq \mathcal F_{\theta, q} \subseteq \mathcal{F}_{\theta, p} $. 
Now, if we let 
\[ \bar{A}_{\theta, p} = \{x\in A_0 + A_1: \exists f\in \mathcal F_{\theta, p} : f(\theta)=x\}\]
endowed with the norm 
\[ \|x\|_{\bar{A}_{\theta, p}} = \inf \{ \|f\|_{\mathcal F_{\theta, p}}: f(\theta)=x \},\] then  $
\bar{A}_{\theta, p} = \bar{A}_\theta$.  The corresponding Schechter spaces $\bar{A}_{\delta'(\theta), p}$ and $\bar{A}^{\delta'(\theta), p}$ are defined in an analogous manner.

\par We now introduce the following auxiliary norms based on \cite{RW}. Given $(u,v)\in \bC \times \bC$, we define
\[ |(u,v)|_{(\bC; \lambda\bC)_\theta}^{(p)} = \inf\{ \|F\|_{\mathcal F_{\theta, p}(\bC, \lambda\bC)} : F(\theta)=u, \, F'(\theta) = \varphi'(\theta)v\}. \]
If $\lambda=1$, we shall simply write $|(u,v)|_\theta^{(p)}$. An adaptation of \cite[\S3A]{RW} yields:

\begin{prop} \label{prop:C-C} $|(u,v)|_\theta^{(\infty)} \leq 1$ if and only if $|u|^2 +|v| \leq 1$. Consequently, 
	\[ |(u,v)|_\theta^{(\infty)}  = \frac{1}2\left( |v| + \sqrt{|v|^2+4|u|^2}\right).\]
\end{prop}
\begin{proof}	
	\par Given $F\in \mathcal F(\bC, \bC) = A(\mathbb S)$ with $\|F\|_{\mathcal F(\bC, \bC)}\leq 1$, the Schwarz-Pick lemma applied to $F\circ \phi$ yields that, for every $z\in \mathbb D$, 
	\[ (1-|z|^2)|F'[\phi(z)]\phi'(z)| \leq 1 - |{F[\phi(z)]}|^2, \]
	and if we substitute $z=0$, 
	\[ |F'(\theta)| \leq |\varphi'(\theta)|(1-|F(\theta)|^2).\]
	Hence, if $F(\theta)=u$ and $F'(\theta)=\varphi'(\theta)v$, then $|u|^2+|v|\leq 1$. Conversely, let $(u,v)\in \bC\times\bC$ with $|u|^2+|v| \leq 1$. Assume first that $|u|^2+|v|=1$, and let
	\begin{equation} \label{eq:psi} \Psi: \mathbb D \to \mathbb D, \quad  \Psi(z) = \begin{cases} \dfrac{ \frac{v}{|v|}z+u}{ 1+\bar u\frac{v}{|v|}z} & \text{ if } v\neq 0 \\ u & \text{ if } v=0, \end{cases} \end{equation}
	which is a holomorphic function satisfying $\Psi(0)=u$ and $\Psi'(0)=v$. Observe that, if $v\neq 0$, $\Psi$ is just the composition of the rotation $z \mapsto \frac{v}{|v|}z$ with a conformal map from $\mathbb D$ into itself that takes $0$ to $u$. Then, $F=\Psi \circ \varphi \in \mathcal F(\bC, \bC)$ satisfies $F(\theta)=u$, $F'(\theta)=\varphi'(\theta)v$ and $\|F\|_{\mathcal F(\bC, \bC)} \leq 1$. Therefore, $\|(u,v)\|_\theta^{(\infty)} \leq 1$. Now, if $|u|^2 +|v|<1$, then there is $\lambda> 1$ such that $|\lambda u|^2 + |\lambda v| =1$. Hence $\|(\lambda u, \lambda v)\|_\theta^{(\infty)} \leq 1$ and so $\|(u,v)\|_\theta^{(\infty)} \leq \lambda^{-1} < 1$. 
\end{proof}

\par Let us now consider the case $\lambda\neq 1$.

\begin{prop} \label{prop:C-lC}
	Let $p\in [1, +\infty)$. Given $(u,v)\in \bC\times \bC$, then 
	\[ |(u,v)|_\subC{p} = \l^\theta |(u,w)|_\theta^{(p)} \]
	where $w=v+\log \l \phi'(0)u$. 
	\par In particular, 
	\[ |(u,v)|_\subC{\infty} = \frac{\lambda^\theta}{2}\left(|{v+\log\lambda\,  \phi'(0)u}| + \sqrt{|v+\log\lambda\, \phi'(0)u|^2 + 4|u|^2}\right)\]
\end{prop}
 
\begin{proof} The map 
	\[ \cF_{\theta, p}(\bC, \l\bC) \to  \cF_{\theta, p}(\bC, \l\bC) \quad, \quad F(z) \mapsto G(z)=\l^z F(z)\]
	defines a linear isometry in such a way that if $F(\theta)=u$ and $F'(\theta)=\varphi'(\theta)v$, then $G(\theta) = \l^\theta u$ and $G'(\theta)=\varphi'(\theta)\l^\theta[v+\log \l\,  \phi'(0) u] = \varphi'(\theta)\l^\theta w$. Hence it follows that
	\begin{align*}
		|(u,v)|_\subC{p} &= \inf\{\|F\|_{\cF_{\theta, p}(\bC; \bC)}: F(\theta)=u, F(\theta)=\varphi'(\theta)v\}= \\ &= \inf\{ \|G\|_{\cF_{\theta, p}(\bC; \l\bC)}: G(\theta) = \l^\theta u, \ G'(\theta)=\varphi'(\theta)\l^\theta w\} = \\
		&= \l^\theta |(u,w)|_\theta^{(p)}.
	\end{align*}
	\par In particular, for the case $p=\infty$, Proposition \ref{prop:C-C}, yields
	\begin{align*} |(u,v)|_\subC{\infty} &= \l^\theta |(u,v)|_\theta^{(\infty)} = \\ &= \frac{\l^\theta}{2}\left(|v+\log\l\,\phi'(0)u| + \sqrt{|v+\log\l\,\phi'(0)u|^2 + 4|u|^2}\right). \qedhere 
	\end{align*}
\end{proof}

\begin{thm} \label{thm:h-H} The following holds:
	\begin{enumerate}
		\item[i)] $(\bC, \l\bC)_{\delta'(\theta)} \equiv h(\l,\theta)\bC$, where
		\[ h(\l, \theta) = \begin{cases}
			\dfrac{\l^\theta}{|{\log\l}|} & \text{ if } |{\log \l}|\geq \dfrac{\pi}{\sin \pi\theta}, \\[3mm]
			\l^\theta\dfrac{2\pi\sin\pi\theta}{\pi^2+|{\log \l}|^2(\sin \pi\theta)^2} & \text{ if }  |{\log \l}|< \dfrac{\pi}{\sin \pi\theta}.
		\end{cases}\]
		
		\item[ii)] $(\bC, \l\bC)^{\delta'(\theta)} \equiv H(\l,\theta)\bC$, where
		\[ H(\l, \theta) = \frac{\l \sin\pi\theta}{\pi}\left( |{\log\l}| + \sqrt{\frac{\pi^2}{\sin^2 \pi\theta} + |{\log \l}|^2} \right).\]
	\end{enumerate} 
\end{thm}
\begin{proof} Applying Proposition \ref{prop:C-lC}, one gets
	\begin{align*}
		\|1\|_{\delta'(\theta)} &= \inf\{\|F\|_{\cF(\bC, \l\bC)}: F'(\theta)=1\} = \\ &= \inf_{u\in \bC} \inf\{\|F\|_{\cF(\bC, \l\bC)} : F(\theta)= u,\, F'(\theta)=1\} \\
		& = |\phi'(0)| \, \inf_{u\in \bC} \inf\{\|F\|_{\cF(\bC, \l\bC)} : F(\theta)= u,\, F'(\theta)=\varphi'(\theta)\} = \\ &= |\phi'(0)|\, \inf_{u\in \bC} \|(u,1)\|_\subC{\infty} = \\
		& = \frac{\l^\theta}{2}|\phi'(0)|\, \inf_{w\in \bC} \left[ |1+w|+\sqrt{|1+w|^2 + \frac{4|\varphi'(\theta)|^2}{|{\log\l}|^2}|w|^2}\, \right].
	\end{align*}
	All that is left is to solve the above extreme problem, which is equivalent to the following (easier) one: 
	\[ \inf_{a\in \mathbb R} \left[|1+a| + \sqrt{(1+a)^2 + Ma^2}\right] = \begin{cases}
		\sqrt{M} & \text{ if } 0 \leq M \leq 1, \\[2mm]
		\dfrac{2M}{M+1} & \text{ if } M>1, 
	\end{cases} \]
	where 
	\[ M=\frac{4|\varphi'(\theta)|^2}{|{\log \l}|^2} = \frac{\pi^2}{(\sin \pi\theta)^2|{\log \l}|^2}.\]
	This suffices to show $(i)$. The proof of $(ii)$ is analogous, using the fact that $\|1\|^{\delta'(\theta)} = \|(1,0)\|_\subC{\infty}$.
\end{proof}

\par We now prove that $\bar A_{\delta(\theta), p} = \bar{A}_{\delta'(\theta)}$ for the simple couple $\bar A=(\bC, \l\bC)$. Actually, this is true for a general couple, but we only need this particular case. Our next result make use of the \emph{Hardy spaces} on the disk $H^P(\mathbb D)$, which consist of holomorphic functions $f: \mathbb D \to \mathbb C$ such that, for $1\leq p< \infty$,
\[ \|f\|_p = \sup_{0\leq r<1} \left( \int_{0}^{2\pi} f(re^{i\theta}) \, \frac{d\theta}{2\pi} \right)^\frac{1}{p}< +\infty,\]
and, for $p=\infty$,
\[ \|f\|_\infty = \sup_{|z|<1} |f(z)| < \infty.\] 

\begin{prop} \label{prop:p} For any $p\geq 1$, 
	\begin{enumerate}
		\item[i)] $(\bC, \l\bC)^{\delta'(\theta), p} = (\bC, \l\bC)^{\delta'(\theta)}$.
		\item[ii)] $(\bC, \l\bC)_{\delta'(\theta), p} = (\bC, \l\bC)_{\delta'(\theta)}$.
	\end{enumerate}
\end{prop}
\begin{proof} Recall that Proposition  \ref{prop:C-lC} provides the identities
	\[	\|1\|_{\delta'(\theta)} = |\phi'(0)|\, \inf_{u\in \bC} |(u,1)|_\subC{\infty}, \qquad
	\|1\|^{\delta'(\theta)} = |(1,0)|_\subC{\infty}, \]
	and for any $(u,v)\in \bC \times \bC$, Proposition \ref{prop:C-C}, in its turn, yields:
	\[ |(u,v)|_\subC{p} = \l^\theta |(u,w)|_{\theta}^{(p)}, \quad w=v +\log\l\cdot \phi'(0)u. \]
	Now, observe that the very definition of the Poisson kernels on the strip implies that $F\in \mathcal F(\bC, \bC)$ precisely when $F \circ \varphi \in H^p(\mathbb D)$, and $\|F\|_{\mathcal F_{p,\theta}} = \|F\circ \varphi \|_{H^p(\mathbb D)}$. This yields
	\[ |(u,v)|_\theta^{(p)} = \inf\{\|G\|_{H^p(\mathbb D)}: G(0)=u, \, G'(0)=v\}. \]
	\par We now combine this information appropriately. Applying the Cauchy integral formula to any $G\in H^1(\mathbb D)$ such that $G(0)=u$ and $G'(0)=v$ we obtain that
	\begin{align*}  |(u,v)|_\subC{\infty} & = \frac{\l^\theta}2\left(|w|+\sqrt{|w|^2+4|u|^2}\right) \lesssim \l^\theta\max\{|u|, |w|\} \leq \l^\theta |(u,w)|_\theta^{(1)} = \\[2mm] &= |(u,v)|_\subC{1}.  \end{align*}
	Finally, the inequalities
	\[ |(u,v)|_\subC{1} \lesssim |(u,v)|_\subC{p} \lesssim  |(u,v)|_\subC{\infty}\] 
	follow easily from the fact that  $\mathcal F \subseteq \mathcal F_{\theta, p} \subseteq \mathcal F_{\theta, 1}$ for every $p>1$. 
\end{proof}
\par We now focus on the case $p=2$, which is particularly simple. 

\begin{thm}  The following holds:
	\begin{enumerate}
		\item[i)] $(\bC, \l\bC)_{\delta'(\theta), 2} \equiv h_2(\l,\theta)\bC$, where
		\[ h_2(\l, \theta) = \l^\theta \frac{\frac{2\sin\pi\theta}{\pi}}{\left(1+\left|\frac{2\sin\pi\theta}{\pi}\log\l\right|^2\right)^\frac12}.\]
		
		\item[ii)] $(\bC, \l\bC)^{\delta'(\theta), 2} \equiv H_2(\l,\theta)\bC$, where
		\[ H_2(\l, \theta) = \l^\theta \left(1+\left|\frac{2\sin\pi\theta}{\pi}\log\l\right|^2\right)^\frac12.\]
	\end{enumerate} 
\end{thm}
\begin{proof}
	It is a simple exercise to see that 
	\[ |(u,v)|_\theta^{(2)} = \sqrt{|u|^2+|v|^2}. \]
	Now we use again the reasoning of Proposition \ref{prop:C-lC} to obtain
	\[ |(u,v)|_\subC{2} = \l^\theta\left(|u|^2 + |v+\log\l\,\phi'(0)u|^2\right)^\frac12, \]
	and then proceed as in the proof of Theorem \ref{thm:h-H} to obtain
	\begin{align*}
		\|1\|_{(\bC, \l\bC)_{\delta'(\theta), 2}} &= |\phi'(0)|\, \inf_{u\in \bC} \l^\theta\left(|u|^2+|{1+\log\l\,\phi'(0)u}|^2\right)^\frac12 = \\[2mm]
		& = \l^\theta \frac{|\phi'(0)|}{\big(1+|{\log\l\,\phi'(0)}|^2\big)^\frac12}.
	\end{align*}
	 ]
	The second isometry follows similarly:
	\[
	\|1\|^{(\bC; \l\bC)_{\delta'(\theta), 2}} = |(1,0)|_\subC{2} = \l^\theta\left(1+|{\log\l\,\phi'(0)}|^2\right)^\frac12. \qedhere \]
\end{proof}

\subsection{The couple $(L^p(w_0), L^p(w_1))$}
\par We finally obtain the desired information regarding the endpoint spaces for the couple $(L^p(w_0), L^p(w_1))$, where $w_1 \in N$. 
The identification of the space  $(L^p(w_0), L^p(w_1))_{\delta(0)}$ bears no difficulty if we recall that for any $\theta \in (0,1)$ we have
\[ (L^p(w_0), L^p(w_1))_\theta = L^p(w_0^{1-\theta}w_1^\theta).\]
Hence it is easy to obtain the following result.
\begin{prop}
	If $D$ is a dense subset of $L^p(w_0)\cap L^p(w_1)$ and $w_1\in N$, then the norms $\|\cdot\|_{\delta(0)}$ and $\|\cdot\|_{L^p(w_0)}$ are equivalent on $D$. In particular, the inclusion ${(L^{p}(w_0), L^{p}(w_1))}_{\delta(0)} \hookrightarrow L^{p_0}$ is bounded. 
\end{prop}

\begin{remark}
	\label{remark:main} In order to proceed with the identification of $(L^p(w_0), L^p(w_1))_{\delta'(\theta)} $ and $(L^p(w_0), L^p(w_1))^{\delta'(\theta)}$ for $\theta\in (0,1)$, we need some precisions concerning the measurability of certain functions. Let us denote $\mathcal M$ the underlying measure space of the spaces $L^p(w_0)$ and $L^p(w_1)$. 
	\begin{enumerate}
		\item Given $u(x)$ and $v(x)$ two measurable functions, there is a function $F_1(z,x)$ such that:
		\begin{itemize}
			\item For every $z\in \mathbb S$, $x\to F_1(z,x)$ is measurable.
			\item For every $x\in \mathcal M$, $z\to F_1(z,x) \in \mathcal F(\bC, \bC)$, with $F_1(\theta, x) = u(x)$, $F_1'(\theta, x) = \varphi'(\theta)\, v(x)$ and $\|F_1(\cdot, x)\|_{\cF(\bC; \bC)} \leq |(u(x), v(x))|_\theta^{(\infty)}$.
		\end{itemize}
		
		\vspace*{1mm}
		
		\par \noindent Indeed, let $\mu(x) = |(u(x), v(x))|_\theta^{(\infty)} = \frac{1}{2}(|v(x)| + \sqrt{|v(x)|^2+4|u(x)|^2})$ and set $U(x) = {u(x)}/{\mu(x)}$, $V(x) = {v(x)}/{\mu(x)}$. Then, $F_1(z,x) = \mu(x)\cdot \Psi(\varphi_\theta(z), x)$, where 
		\[ \Psi: \mathbb D \times \mathcal M \to \mathbb C, \quad \Psi(\xi, x)  = \begin{cases} \dfrac{ \frac{V(x)}{|V(x)|}\xi+U(x)}{ 1+\overline {U(x)}\frac{V(x)}{|V(x)|}\xi} & \text{if } V(x)\neq 0 \\[4mm] 0 & \text{if } V(x)=0, 
			\end{cases}
		\]
		is the corresponding map given in \eqref{eq:psi}.
		
		\item Then, since in Proposition \ref{prop:C-lC} it is proved that, for any $u,v\in \bC$,  
		\[ |(u,v)|_\subC{\infty} = \l^\theta |(u,w)|_\theta^{(\infty)}\]
		where $w=v+\phi'(0)\log\l \, u$, we conclude that if $u(x)$, $v(x)$ and $\l(x)$ are three measurable functions, there exists a function $F_2(z,x)$ such that 
		\begin{itemize}
			\item For every $z\in \mathbb S$, $x\to F_2(z,x)$ is measurable.
			\item For every $x\in \mathcal M$, $z\to F_2(z,x) \in \mathcal F(\bC; \l(x) \bC)$, with $F_2(\theta, x) = u(x)$, $F_2'(\theta, x) = \varphi'(\theta) v(x)$ and such that $\|F_2(\cdot, x)\|_{\cF(\bC; \l(x)\bC)} \leq |u(x), v(x)|_\subC{\infty}$.
		\end{itemize}
		We just need to apply the previous remark to the functions $\tilde u(x) = \l(x)^\theta u(x)$ and $\tilde v(x) =\l(x)^\theta\big(v(x) + \varphi'(0)\log\l(x)\, u(x)\big)$ to obtain a function $\widetilde F_1(z,x)$ and then consider $F_2(z,x) = \l(x)^{-z}\, \widetilde F_1(z,x)$.
		
		\item Finally, in Theorem \ref{thm:h-H} it is shown that if $\l(x)$ is measurable, then 
		\[ h(\l(x), \theta) = \|1\|_{{(\bC; \l(x)\bC)}_{\delta'(\theta)}} \geq \phi'(0) |(u(x), 1)|_{{(\bC; \l(x)\bC)}_\theta}^{(\infty)} \]
		for every measurable function $u(x)$. 
		Therefore, there exists a function $F(z,x)$ such that 
		\begin{itemize}
			\item For every $z\in \mathbb S$, $x\to F(z,x)$ is measurable.
			\item For every $x\in \mathcal M$, $z\to F(z,x) \in \mathcal F(\bC, \l(x)\bC)$, with $F'(\theta, x) = 1$ and $\|F(\cdot, x)\|_{\cF(\bC; \l(x)\bC)} \leq h(\l(x), \theta)$.
		\end{itemize}
		By the same token, using that 
		\[ H(\lambda(x), \theta) = \|1\|^{\delta'(\theta)} = \|(1,0)\|_\subC{\infty},\]
		we infer the existence of a function $F(z,x)$ such that
		\begin{itemize}
			\item For every $z\in \mathbb S$, $x\to F(z,x)$ is measurable.
			\item For every $x\in \mathcal M$, $z\to F(z,x) \in \mathcal F(\bC, \l(x)\bC)$, with $F(\theta, x) = 1$, $F'(\theta, x) = 0$ and $\|F(\cdot, x)\|_{\cF(\bC; \l(x)\bC)} \leq H(\l(x), \theta)$.
		\end{itemize}

	\end{enumerate}
\end{remark}

\begin{thm} Let $p\in [1, +\infty)$. Then
	\begin{itemize}
		\item $(L^p(w_0), L^p(w_1))_{\delta'(\theta)} = L^p\left(w_0\, h\left(\frac{w_1^{1/p}}{w_0^{1/p}}, \theta\right)^p\right)$. \vspace*{2mm}
		\item $(L^p(w_0), L^p(w_1))^{\delta'(\theta)} = L^p\left(w_0\, H\left(\frac{w_1^{1/p}}{w_0^{1/p}}, \theta\right)^p\right)$.\vspace*{2mm}
	\end{itemize}
	\par In particular, for $p=2$,
	\begin{itemize}
		\item $(L^2(w_0), L^2(w_1))_{\delta'(\theta),2} \equiv L^2\left(w_0\, h_2\left(\frac{w_1^{1/2}}{w_0^{1/2}}, \theta\right)^2\right)$. \vspace*{2mm}
		\item $(L^2(w_0), L^2(w_1))^{\delta'(\theta),2} \equiv L^2\left(w_0\, H_2\left(\frac{w_1^{1/2}}{w_0^{1/2}}, \theta\right)^2\right)$.
	\end{itemize}
\end{thm}
\begin{proof} During the proof, we shall write $\l(x) = \frac{w_1(x)^{1/p}}{w_0(x)^{1/p}}$, and $\l(x)=0$ if $w_0(x)=0$. 
	Let us prove the first equality. By Remark \ref{remark:main}, there is a function $F\in \cF(\bC; \l(x)\C)$ such that $F(z, \cdot)$ is measurable for every $z\in \mathbb S$, $F'(\theta, x) = 1$, and $\|F\|_{\cF(\bC; \l(x)\bC)} \leq h(\l(x), \theta)$. Therefore, given $f\in L^p(w_0\, h(\l, \theta)^p))$, if one considers $G(z,x) = f(x) F(z,x)$, one can easily see that $G\in \cF(L^p(w_0), L^p(w_1))$ with $\|G\|_{\cF(L^p(w_0), L^p(w_1))} \leq \|f\|_{L^p(w_0\, h(\l, \theta)^p)}$, and that $G'(\theta, x) = f(x)$. Hence $f\in {(L^p(w_0), L^p(w_1))}_{\delta'(\theta)}$ and 
	\[ \|f\|_{(L^p(w_0), L^p(w_1))_{\delta'(\theta)}} \leq \|f\|_{L^p(w_0\, h(\l, \theta)^p)}. \]
	For the converse, observe that if $F\in \cF(L^p(w_0), L^p(w_1))$, then an application of Fubini's theorem yields
	\[ \|F\|_{\cF_{p, \theta}} = \left( \int_{\mathcal M} \|F(\cdot, x) \|^p_{\cF(\bC; \l(x)\bC)_{\theta, p}} w_0(x) \, dx \right)^\frac1p. \]
	Hence,
	\begin{align*}
		\|f\|_{(L^p(w_0), L^p(w_1))_{\delta'(\theta)}} &\geq \|f\|_{(L^p(w_0), L^p(w_1))_{\delta'(\theta), p}} = \\[1mm] &
		= \inf_{F'(\theta, x) = f(x)} \left( \int_{\mathcal M} \|F(\cdot, x)\|^p_{\cF(\bC; \l(x)\bC)_{\theta, p}} w_0(x) \, dx \right)^\frac1p \geq  \\[1mm]
		& \geq \left( \int_{\mathcal M}  \inf_{F'(\theta, x) = f(x)} \|F(\cdot, x)\|^p_{\cF(\bC; \l(x)\bC)_{\theta, p}} w_0(x) \, dx \right)^\frac1p = \\[1mm] 
		& = \left( \int_{\mathcal M} \|f(x)\|^p_{{(\bC; \l(x)\bC)}_{\delta'(\theta), p}} w_0(x) \, dx \right)^\frac1p \gtrsim\\[1mm] 
		& \overset{(\star)}{\gtrsim} 
		\left( \int_{\mathcal M} |f(x)|^p w_0(x) h(\l(x), \theta)^p \, dx\right)^\frac1p  		= \|f\|_{L^p(w_0\, h(\l, \theta)^p)}.
	\end{align*}
	where the step marked with $(\star)$ holds by virtue of Proposition \ref{prop:p} and Theorem \ref{thm:h-H}. 
	The remaining identities are proved in a similar fashion. 
\end{proof}

\par We finally arrive to the desired result concerning the endpoint spaces: 
\begin{thm} Let $p\in [1, +\infty)$ and $D$ be a dense subset of $L^p(w_0) \cap L^p(w_1)$, where $w_1 \in N$. Then 
	\begin{enumerate}
		\item[i)] The norms $\|\cdot\|^{\delta'(0)}$ and $\|\cdot\|_{L^p(w_0)}$ are equivalent on $D$. 
		\item[ii)] The natural inclusion $(L^p(w_0), L^p(w_1))_{\delta'(0)} \hookrightarrow L^p(w_0)$ is bounded. 
	\end{enumerate}
\end{thm}
\begin{proof} Let us write once again $\l(x) = \frac{w_1(x)^{1/p}}{w_0(x)^{1/p}}$ and $\l(x)=0$ if $w_0(x)=0$. To check ($i$), it is enough to observe that, since $w_1(x)\neq 0$ a.e.$(x)$, then 
	\[ \lim_{\theta\to 0^+} w_0 \, H(\l, \theta)^p = w_0 \quad \text{ a.e.} (x).\]
	Hence, if $d\in D$, the previous theorem asserts that
	\[ \|d\|^{\delta'(0)} = \limsup_{\theta \to 0^+} \|d\|^{\delta'(\theta)} \approx \lim_{\theta \to 0^+} \|d\|_{L^p(w_0H(\l, \theta)^p)} = \|d\|_{L^p(w_0)}. \]
 	
	\par The proof of ($ii$) follows an analogous argument. First, appealing again to the fact that $w_1(x)\neq 0$ a.e$(x)$, we obtain that 
	\[ \lim_{\theta\to 0^+} w_0\, \frac{h(\lambda, \theta)}{\theta} \approx w_0 \quad \text{a.e.} (x). \]
	Now, given $f\in (L^p(w_0), L^p(w_1))_{\delta'(0)}$, there exists a decreasing sequence $(\theta_k)_{k=1}^\infty$ converging to $0$ in such a way that $f\in \bigcap_{k=1}^\infty (L^p(w_0), L^p(w_1))_{\delta'(\theta_k)}$ and 
	\[ \lim_k \frac{\|f\|_{(L^p(w_0), L^p(w_1))_{\delta'(\theta_k)}}}{\theta_k} = \|f\|_{(L^p(w_0), L^p(w_1))_{\delta'(0)}}. \]
	This implies that
	\begin{align*}
		\|f\|_{L^p(w_0)} & =  \left(\int_{\mathcal M} |f(x)|^pw_0(x)\, dx\right)^\frac1p \lesssim \\ & \lesssim \left( \int_{\mathcal M} \liminf_k \ \frac{1}{\theta_k} |f(x)|^pw_0(x) h(\l(x), \theta_k)^p\, dx \right)^\frac1p \leq \\[1mm]
		& \leq \left( \liminf_k \int_{\mathcal M}  \ \frac{1}{\theta_k} |f(x)|^pw_0(x) h(\l(x), \theta_k)^p\, dx \right)^\frac1p \leq \\[1mm]
		& \leq \left( \limsup_k \int_{\mathcal M}  \ \frac{1}{\theta_k} |f(x)|^pw_0(x) h(\l(x), \theta_k)^p\, dx \right)^\frac1p =  \\[1mm]
		&=\limsup_k \frac{1}{\theta_k} \|f\|_{L^p(w_0 h(\l, \theta_k)^p)} \lesssim \\[1mm]
		& \lesssim \|f\|_{(L^p(w_0), L^p(w_1))_{\delta'(0)}}, 
	\end{align*}
	where the last step is guaranteed by our previous theorem. This is enough to conclude.
  
\end{proof}

\begin{remark} The arguments given in this section can be improved to show that, for the couples $\bar{A} = (L^{p_0}, L^{p_1})$ and $\bar{A} = (L^p(w_0), L^p(w_1))$ for $1\leq p_0, p_1, p \leq \infty$ and $w_1 \neq 0$ a.e., the completion of \emph{all} the endpoint spaces coincides with the base space $A_0$. However, observe that this does not improve the situation regarding the applications of the paper, since for a function $m$ to be in $\mathcal M_p$, we require that $T_m$ maps $L^p$ to $L^p$, and not a space isomorphic to $L^p$ to a space isomorphic to $L^p$. 
\par Therefore, we believe that a deeper study of the endpoint spaces $\bar{A}^{\delta'(0)}$ and $\bar{A}_{\delta'(0)}$ for a general interpolation couple $\bar{A}$ is in order. In particular, it seems reasonable to ask under which hypotheses it can be ensured that the completion of the endpoint spaces is isomorphic to $A_0$.
\end{remark}

\end{document}